\documentclass{article}

\usepackage{bm}
\usepackage{graphicx}
\usepackage{tabulary,booktabs}

\newcommand{\libmesh}[0]{{\tt libMesh}}
\newcommand{\petsc}[0]{{\tt PETSc}}

\newcommand{\radiusamr}{r_{\mbox{\tiny AMR}}}
\newcommand{\rmin}{r_{\mbox{\tiny min}}}

\newcommand{\W}{\Omega}

\newcommand{\colvec}[2]{\left(\!\!\!\begin{array}{c} #1 \\ #2 \end{array}\!\!\!\right)}

\newcommand{\matt}[4]{\left[\begin{array}{cc} #1 & #2 \\ #3 & #4 \end{array}\right]}

\newcommand{\comment}[1]{} 

\begin{document}

\title{Dynamic Adaptive Mesh Refinement for Topology
Optimization\thanks{This work was supported in part by the
National Science Foundation under Grant DMR-03 25939 ITR
through the Materials Computation Center at the University of
Illinois at Urbana-Champaign.}}

\author{
Shun Wang\thanks{Department of Computer Science, University of Illinois at
Urbana-Champaign, Urbana, Illinois 61801, U.S.A., wangshun98@gmail.com,},
Eric de Sturler\thanks{Department of Mathematics, Virginia Tech,
Blacksburg, Virginia 24061, U.S.A., sturler@vt.edu,},
Glaucio H. Paulino\thanks{Department of Civil and Environmental Engineering,
University of Illinois at Urbana-Champaign, Urbana, Illinois 61801, U.S.A., paulino@uiuc.edu.}
}
\date{}
\maketitle

\begin{abstract}
We present an improved method for topology optimization with
both adaptive mesh refinement and derefinement. Since the total
volume fraction in topology optimization is usually modest,
after a few initial iterations the domain of computation is
largely void. Hence, it is inefficient to have many small
elements, in such regions, that contribute significantly to the
overall computational cost but contribute little to the
accuracy of computation and design. At the same time, we want
high spatial resolution for accurate three-dimensional designs
to avoid postprocessing or interpretation as much as possible.
Dynamic adaptive mesh refinement (AMR) offers the possibility
to balance these two requirements. We discuss requirements on
AMR for topology optimization and the algorithmic features to
implement them. The numerical design problems demonstrate (1)
that our AMR strategy for topology optimization leads to
designs that are equivalent to optimal designs on uniform
meshes, (2) how AMR strategies that do not satisfy the
postulated requirements may lead to suboptimal designs, and (3)
that our AMR strategy significantly reduces the time to compute
optimal designs.
\end{abstract}

\vspace{11pt}
\noindent
{\bf keywords:} adaptive mesh refinement, topology optimization, iterative solvers.

\section{Introduction}

Topology optimization is a powerful structural optimization method that combines a numerical solution method, usually the
finite element method, with an optimization algorithm to find the optimal material distribution inside a given domain
\cite{RahmSwan2004, Sigmund2000, PaulinoPS2005, PaulLe2008}. In designing the topology of a
structure we determine which points in the domain should be material and which points should be void. However, it is
known that, in the continuum setting, topology optimization leads to designs with intermediate densities. So
continuous values between 0 and 1 replace discrete values (0 or 1) to represent the relative densities, and some form of
penalization is used to obtain designs with almost discrete 0/1 material density distribution \cite{Bendsoe1989SIMP}.

In topology optimization, problems are solved most commonly on
fixed uniform meshes with a relatively large number of elements
in order to achieve accurate designs \cite{Rozvany2001,
Mackerle2003}. However, as void and solid regions appear in the
design, it is more efficient to represent the holes with fewer
large elements and the solid regions, especially the material
surface, with more fine elements. Since the shape and position
of holes and solid regions are initially unknown, the most
economical mesh representation for the design is unknown {\it a
priori}. Therefore, adaptive mesh refinement (AMR) is very
suitable for topology optimization. {\em The purpose of AMR for
topology optimization is to get the design that would be
obtained on a uniformly fine mesh, but at a much lower
computational cost by reducing the total number of elements and
having fine elements only where and when necessary}.

Highly accurate designs on uniform meshes may require so many
elements that the solution of the optimization problem becomes
intractable. However, AMR leads to high resolution in the mesh
only when and where necessary. This makes it possible to obtain
accurate designs with a modest number of elements and hence
with a reasonable cost. Even when a design on a uniform mesh is
computationally feasible, AMR tends to reduce the computational
cost by reducing the required number of elements and by
improving the conditioning of linear systems arising from the
finite element discretization. Obviously, we do not want the
use of AMR or the AMR procedure to alter the computed designs.
However, there is a risk of this, since the mesh influences the
computed deformations and sensitivities. It is imperative then
that the solutions from the finite element analysis using AMR
must be as accurate as those obtained on a uniform fine mesh.
\footnotemark \footnotetext{For comparison everywhere in this
paper the element size for the uniform fine mesh is the same as
the element size at the highest level of refinement in the AMR
mesh.}
Moreover, the final design must be governed by accurate
sensitivities corresponding to those obtained on the finest
mesh. If coarse mesh solutions drive or limit the design,
suboptimal designs may result when designs optimal on a coarser
mesh differ substantially from the optimal design on a (much)
finer mesh. We will demonstrate that this occurs in quite
simple cases. The early work in this area, though leading to
acceptable designs in specific instances, does not satisfy
these properties. We will propose relatively simple but
essential changes to these methodologies that lead to AMR-based
designs that are equivalent (up to some small tolerance) to
designs on uniform fine meshes. In addition, our approach leads
(in principle) to an overall more efficient method as we reduce
the total number of elements further. The topology optimization
may lead to a sequence of (intermediate) structures requiring
high mesh resolution in different parts of the computational
domain. Therefore, it is important to (1) allow the meshes at
all levels to change continually (dynamic) and (2) to allow
both mesh refinement and derefinement \cite{Wang2007}.
Derefinement is important for efficiency when the initial
discretization needs to include relatively small elements in
certain regions. This is important in a number of cases, which
are elaborated upon below.

In the next section, we provide an assessment of previous AMR
strategies, namely the implementations by Costa and Alves
\cite{Costa2003} and Stainko \cite{stainko2006}. In
Section~\ref{sec:top_opt}, we provide a brief introduction to
topology optimization. Next, in Section~\ref{sec:dynamic_amr},
we state the purpose of our AMR strategy for topology
optimization and explain the requirements it poses. Based on
these requirements, we propose a more robust and dynamic AMR
strategy. We describe a number of implementation issues of our
AMR strategy in Section~\ref{sec:impls}. We briefly discuss the
iterative solution of the large sparse linear systems arising
in the finite element analysis in Section~\ref{sec:iter_sol}.
In Section~\ref{sec:results}, we show numerical experiments
that demonstrate the robustness and efficiency of our AMR
strategy. The first experiment also explains why the refinement
strategies by Costa and Alves \cite{Costa2003} and Stainko
\cite{stainko2006} may lead to suboptimal designs. Finally, in
Section~\ref{sec:conc}, we present conclusions about our AMR
strategy for topology optimization algorithms, and we mention
some directions for future work.

\section{Assessment of Previous AMR Strategies} \label{sec:prev_amr}

Little research has been done in applying AMR to topology
optimization. So, we start by briefly discussing two recent,
important, papers in this area. The AMR method by Costa and
Alves \cite{Costa2003} goes through a predetermined, fixed
sequence of optimizations and subsequent mesh refinements (they
do not use derefinements), always using (or assuming) a
converged solution on a `coarse mesh' to guide the refinement
of that mesh and start the optimization on the next `fine
mesh'. Coarse meshes and the solutions on these coarse meshes
are never revisited or updated after generating the next finer
mesh. The method aims at refining the coarse mesh design. Hence
the region with the fine(r) mesh that contains the material
boundary is always confined to the coarse mesh region that has
been determined before using only coarse mesh calculations.
After a fixed number of optimization steps on a given mesh,
they refine all material elements (density $0.5$ or larger) and
elements on the boundary between material elements and void
elements (density less than $0.5$). Furthermore, they refine
elements that do not satisfy certain quality or solution error
criteria. In addition, there are a few special cases that also
lead to refinement. These refinements lead to accurate finite
element solutions in material regions, a high mesh resolution
on the material boundary and, therefore, accurate
representation of this boundary, and larger elements in void
regions reducing the amount of work. However, as reported by
the authors, the `optimal design' found by the method depends
on the initial mesh and is not the same as the optimal design
found using a uniform fine mesh \cite{Costa2003}. Although the
authors do not report this, we conjecture that the design found
using the adaptively refined mesh is not an `optimal design' on
the uniform mesh, that is, it has higher compliance than the
solution obtained on the uniform fine mesh. See also our
numerical experiments below. Finally, only two-dimensional
designs are treated, but conceptually we expect their algorithm
to work similarly in three-dimensional designs.

Stainko follows a slightly different approach with respect to
the refinements \cite{stainko2006}. Mesh refinement is done
only along the material boundary as indicated by the
(regularization) filter. So, elements completely inside a
material region or a void region are not refined. In principle
this leads to a smaller number of elements and hence a reduced
computation time. However, Stainko's procedure also progresses
strictly from coarser meshes to finer meshes, and a coarse mesh
is never updated after mesh refinement. So, just as in
\cite{Costa2003}, the finest mesh, which contains the material
boundary, is always confined to regions defined by (all)
earlier refinements (all refinements are nested), each of which
is based only on the corresponding coarse mesh computations.
Stainko does not test whether the designs obtained are the same
as those obtained on uniformly fine meshes; however, our
experiments below show that, again, the designs will be depend
on the initial mesh (resolution) and are not the same as
optimal designs on the uniformly fine mesh at the maximum
refinement level.

These approaches share two important choices that may lead to
problems. First, both approaches solve the design problem on a
fixed mesh until convergence before carrying out mesh
refinement. After refinement on a given level, the mesh on that
level remains fixed for the remainder of the optimization, and
all further refinements are therefore constrained by the
converged coarser level solutions. This works well in terms of
refining the design, but for many design problems the optimal
solution on a uniform fine(st) mesh is quite different from the
converged solution on a coarser mesh. In that case, mesh
refinement based only on the coarser level solution will
erroneously confine the solution on the finer mesh to a smooth
version of the coarser level solution. Therefore, the
approaches proposed in \cite{stainko2006,Costa2003} may lead to
suboptimal designs, as we will show in our numerical
experiments.

Second, both approaches use only refinement but no
derefinement, which may lead to inefficiencies. First, for
designs with thin structures, the initial, coarsest, mesh must
be fine enough to give a reasonable result. If fine elements
that are no longer required cannot be removed as the design
evolves, then more computational work than necessary will be
performed. Second, in topology optimization approaches that use
filtering for regularization, for an accurate design requiring
a high resolution mesh, the (appropriate) filter will not work
on the coarser meshes, because the filter radius, which should
be a physical feature size independent of the mesh, will
typically be too small. Hence, we must start with a relatively
fine mesh. However, after a modest number of optimization
steps,  large regions will likely have become void and fine
elements could be removed without problems. Again substantial
computational overhead results from having to work with too
many fine elements. Third, any AMR strategy that allows changes
in the design beyond previously computed coarse level designs
and refines the mesh to accommodate such changes will be
inefficient if fine elements in void regions cannot be removed.

Therefore, a more robust and efficient refinement strategy is
needed. Hence, we propose a dynamic meshing strategy that
includes both mesh refinement and derefinement everywhere in
the computational domain. Our improved AMR strategy has two
main components. First, we extend the refinement criteria from
\cite{Costa2003}, refining all material elements and elements
on the boundary, but with an additional layer of refinements
around the material boundary (in the void-region). The
thickness of the layer is a parameter. This way the fine level
design can change shape arbitrarily in optimization steps
between mesh refinements. Second, our AMR method updates coarse
and fine meshes continually, so that small changes arising in
the more accurate computations on finer meshes can change the
shape of the design arbitrarily in the course of the
optimization; the fine(r) meshes move with the material
boundary. This means that our designs are really based on the
accurate fine mesh computations and are not confined to regions
fixed by earlier coarse mesh computations. Since we do
continual mesh adaptation, we may have fine elements in regions
that have become void at some point. Derefinement will remove
those fine elements. Further details are given in the next
subsection. This approach leads to designs on AMR meshes at
greatly reduced cost that are the same as the designs that
would have been obtained on a uniform fine mesh of the highest
resolution (within a small tolerance). We will demonstrate this
experimentally in section~\ref{sec:results}. Our approach also
allows us to start with coarser meshes, since the coarse mesh
solution does not need to be a good approximation to the final
solution. Even when we start with a finer mesh for faster
convergence or proper functioning of the regularization filter,
derefinement allows us to remove fine elements that have become
void (which tends to happen quite quickly).

\section{A Brief Topology Optimization Review} \label{sec:top_opt}

\begin{figure}
\begin{center}
\includegraphics[scale=0.5]{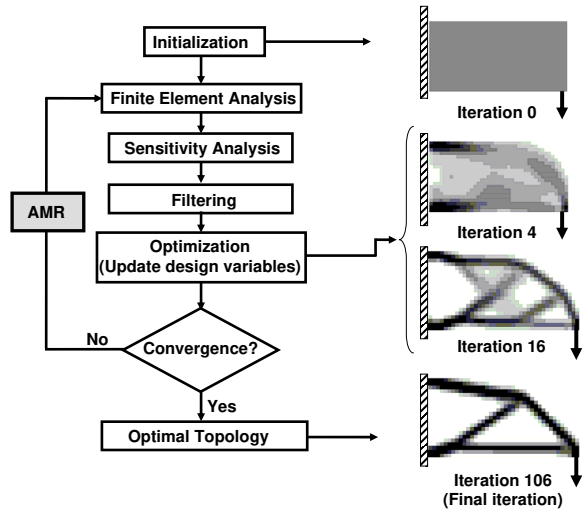}
\end{center}
\caption{Overview of the topology optimization algorithm with dynamic AMR.}
\label{fig:TopOptAlg}
\end{figure}

In topology optimization we solve for the material distribution
in a given computational design domain $\W$. The topology
optimization problem we consider here is to minimize the
compliance of a structure under given loads as a function of
the material distribution. To solve this problem numerically,
we discretize the computational domain using finite elements,
where we usually use a lower order interpolation for the
density field (material distribution) than for the displacement
field. The most common approach (also employed for this paper)
is to use (bi-,tri-)linear interpolation for the displacement
field and constant density in each element. The compliance
minimization problem after finite element discretization is
defined as
\begin{eqnarray}
&& \min_{\rho_e\in[\rho_{\!\!_o},1], \forall e} \bm{f}^T\bm{u} \label{eq:simp} \\
\mbox{s.t.} &&
\left\{ \begin{array}{ll}
\bm{K}(\bm{\rho})\bm{u}=\bm{f} & \mbox{ for } \bm{x} \in \Omega\setminus\Omega_0,\\
\bm{u}=\bm{u}_0 & \mbox{ for } \bm{x} \in \Omega_0, \\
\sum_e \rho_e V_e \leq V_0,
\end{array} \right. \nonumber
\end{eqnarray}
where $\rho_e$ is the density in element $e$, $\bm{\rho}$ is
the vector of element densities, $\bm{K}$ is the stiffness
matrix, a function of the discretized density field
($\bm{\rho}$), $V_e$ is the volume of element $e$, $V_0$ is a
maximum volume (fraction) allowed for the design, and $\W_0$ is
the part of the domain where the displacement is prescribed. To
avoid singularity of the stiffness matrix, we enforce a small
positive lower bound $\rho_{\!\!_o}$ on the element density,
typically $10^{-3}$.

As mentioned in the introduction, our discrete model must drive
the continuous material distribution as much as possible to a
0/1 solution. We use the Solid Isotropic Material with
Penalization (SIMP) method to make the undesirable intermediate
densities between $\rho_{\!\!_o}$ (replacing $0$) and $1$
unfavorable \cite{BendsoeBk2003}. In this case, the elasticity
tensor is defined as a function of the element density,
\begin{equation} \label{eq:stiffness}
\bm{E}_e=\rho_e^p\bm{E}_0,
\end{equation}
where $p$ is the penalization parameter. With $p>1$,
intermediate densities are unfavorable as they provide
relatively little stiffness compared with their material cost.
A common choice is $p = 3$, which results in intermediate
material properties that satisfy the Hashin--Shtrikman bound
for composite materials \cite{HashShtrik63}. To avoid problems
with local minima, we usually apply continuation on the
parameter $p$, that is, we start with $p=1$ and slowly increase
$p$ as the design converges.

The general scheme for topology optimization using AMR is
illustrated in Figure \ref{fig:TopOptAlg}. First, we set up the
geometry, the finite element (FE) mesh, the loading and
boundary conditions, and we initialize the density distribution
$\rho$. Then, we start the optimization loop. In this loop, we
assemble and solve the equilibrium equations
$\bm{K}(\bm{\rho})\bm{u}=\bm{f}$ in (\ref{eq:simp}) using the
FE discretization and a linear solver. Next, in the sensitivity
analysis, we compute the derivatives of the objective function
with respective to the design variables, $\partial c/\partial
\rho_e$. After this, we can apply an optional low-pass filter
to remedy the checkerboard problem \cite{SigmundThesis1994,
Sigmund1997, Sigmund1998}, which can be also addressed by an
alternative minimum length scale approach \cite{Guest2004}. In
the next step, we compute an update of the design variables.
There are various optimization algorithms applicable to
topology optimization. For this paper, we use Optimality
Criteria (OC), a simple approach based on a set of intuitive
criteria \cite{BendsoeBk2003,Bendsoe1988}. After updating the
design variables using a chosen optimization algorithm, we
check the convergence of the design. Under certain conditions,
to be discussed next, dynamic mesh adaptation is carried out
before the (next) finite element analysis.

\section{A Dynamic AMR Strategy}
\label{sec:dynamic_amr}

We base our algorithmic choices on a set of requirements on AMR
codes for topology optimization. As stated above, the purpose
of AMR for topology optimization is to get the design that
would be obtained on a uniform fine mesh, but at a much lower
computational cost by reducing the total number of elements and
having fine elements only where (and when) necessary.

First, since the finite element analysis and the computation of
sensitivities drive the changes in material distribution, they
should be as accurate as on the uniform fine mesh. Therefore,
we need a fine mesh that covers at least the material region
and the boundary. Since the void regions have negligible
stiffness they do not influence the (intermediate) linear
finite element solutions and resulting sensitivity
computations. Thus we do not need a fine mesh inside the void
region, and we can use a refinement criterion similar to that
of Costa and Alves \cite{Costa2003}. At this point we focus on
refinement and derefinement for shape only. Therefore, we are
conservative with respect to accuracy, and we expect that, in
future implementations, good error indicators will lead to
further efficiency gains, in particular because of derefinement
in solid material regions.

Second, the accurate computations on the finest level should
drive the changes in the material distribution. This requires
continual mesh adaptation so that computational results after
refinements can drive updates to the material distribution, and
designs are not confined by earlier coarse grid results. This
also means that as the material region moves close to the
boundary between fine and coarse(r) mesh, additional
refinements allow for further evolution.

Third, we need to ensure that the design can change
sufficiently in between mesh updates. Therefore, we maintain a
layer of additional refinements around the material region (in
the void region) and carry out continual mesh adaptation. Due
to the additional layer of refinements and continual mesh
updates, the design can change arbitrarily following the fine
grid computations and resulting sensitivities, and it is not
confined by earlier coarse grid results. To ensure that the
design accurately reflects the fine mesh computations, we allow
rapid refinements of the mesh early on when voids and material
regions (and hence the boundary) develop, rather than delay
refinements until later stages, when a suboptimal design might
have developed.

Fourth, since the design can change substantially from its
estimate on a coarse mesh, we may have fine elements in void
regions. Those elements must be removed for efficiency,
requiring derefinements of the mesh. To facilitate our strategy
of continual mesh refinement and derefinement, we use a
hierarchical representation of adaptive meshes.

We will now state our refinement and derefinement strategy in more detail.
We adapt the mesh when\\ CASE (i):
\begin{enumerate}
\item the relative change in the compliance is smaller than a given threshold, {\bf and }
\item a given minimum number of optimization steps have been carried out since the last mesh update,
\end{enumerate}
{\bf or} when\\ CASE (ii):
\begin{enumerate}
\setcounter{enumi}{2}
\item a given maximum number of optimization steps have been carried out without meeting conditions 1 and 2.
\end{enumerate}
Condition~1 corresponds to a common convergence criterion for
topology optimization that the maximum change in the design
variables is smaller than a certain tolerance (which we usually
set to 0.01). This condition is satisfied  when the solution is
near a local minimum, which might be caused by a no longer
appropriate mesh. In that case, we must adapt the mesh to allow
the design to change further. If the local minimum is not an
artifact of the mesh, the design will remain the same after
mesh adaptation. Condition~2 prevents refinement and
derefinement from happening too frequently. This is important
as the solution needs to adapt to the changed mesh, so that the
computed sensitivities reflect the design and not mesh
artifacts. This also limits the cost of mesh adaptation. In our
experiments this minimum number of optimization steps is set to
five based on experience, and other values are possible.
Regarding condition~3, in our experiments we adapt the mesh at
least every ten optimization steps. This condition leads to
faster convergence because it ensures that the mesh is
appropriate for the material distribution. Using these
conditions, we can start with a fairly coarse mesh, and we may
carry out mesh (de)refinement before the design converges on
any mesh if necessary.

\vspace{10pt}

\begin{figure}
\begin{center}
\begin{tabular}{cp{0.2in}c}
\includegraphics[scale=0.8]{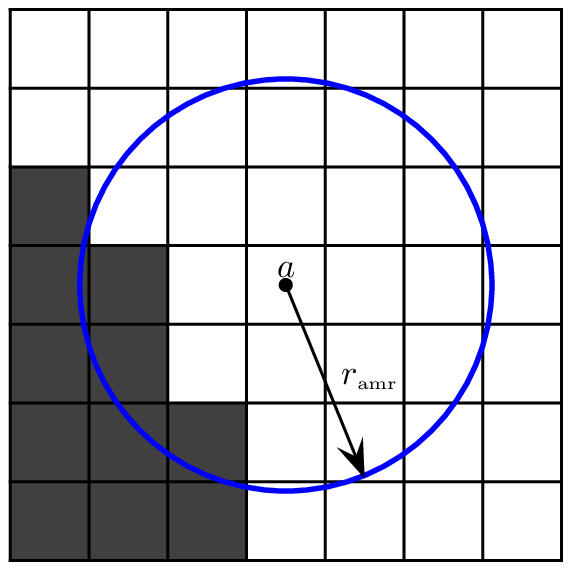} &&
\includegraphics[scale=0.8]{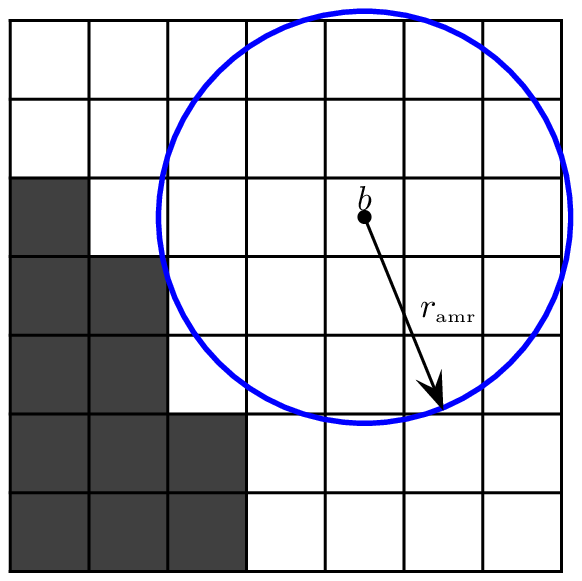}
\end{tabular}
\end{center}
\caption{Refinement criteria for void element. Element $a$ is marked for refinement, because it has
solid elements within distance $\radiusamr$; element $b$ is marked for derefinement.}
\label{fig:amr_radius}
\end{figure}

We adapt our mesh according to the following procedure.
\begin{enumerate}
\item Mark all the elements for refinement or derefinement based on the following
criteria:
\begin{itemize}
\item If element $e$ is solid, i.e., $\rho_e\in[\rho_s,1]$, where $\rho_s$ is a chosen density threshold, or element $e$ is within a
    given distance $\radiusamr$ from a solid element we mark it for refinement.
\item If element $e$ is void, i.e., $\rho_e\in[\rho_{\!\!_o},\rho_s]$, and there are no solid elements within distance $\radiusamr$, we mark element $e$ for derefinement. See Figure \ref{fig:amr_radius}.
\end{itemize}
\item Check compatibility for the mesh that will be generated and make the following adjustments in two sweeps over all elements:
\begin{itemize}
\item In the first sweep, we unmark elements marked for derefinement, if they have a sibling (an element generated by the same  refinement) that is not marked for derefinement.
\item In the second sweep, we unmark elements marked for derefinement, if derefinement would lead to level two or higher edge incompatibility. We allow only level one incompatibility; see Figure \ref{fig:mesh_incomp} and the discussion in Section~\ref{sec:impls}.
\end{itemize}
\end{enumerate}

The above refinement criteria result in a layer of fine
elements on the void side of the solid/void interface that
allows the material to be redistributed locally. If a material
boundary moves near the fine/coarse element interface, mesh
refinement creates a new layer of fine elements around the
current material surface to allow further local redistribution
of the material. On the other hand, if some fine elements
become void, these fine elements are removed by derefinement to
keep the optimization efficient.

\section{Iterative Solution Scheme} \label{sec:iter_sol}

Although AMR significantly reduces the number of DOFs in the
finite element simulation, we still have to solve a sequence of
large linear systems, especially for three-dimensional designs.
Moreover, because of the large difference in density between
distinct parts of the computational domain and the elasticity
tensor given by (\ref{eq:stiffness}), with $p=3$ toward the end
of the optimization, the linear systems are very
ill-conditioned. Hence, proper preconditioning is essential. In
\cite{TopKrylov2006}, we showed how to precondition the linear
systems arising in topology optimization, and we also used {\em
Krylov subspace recycling}\/ \cite{Parks2006} to reduce the
number of iterations over multiple linear systems. We briefly
mention the main ideas here.

To remedy the serious ill-conditioning in topology optimization
problems, we explicitly rescale each stiffness matrix such that
the diagonal coefficients are all equal, as is the case for a
problem with homogeneous density. We rescale the stiffness
matrix $\bm{K}$ by multiplying with a diagonal matrix on both
sides,
\begin{eqnarray}
\nonumber 
  \tilde{\bm{K}} &=& \bm{D}^{-1/2}\bm{K}\bm{D}^{-1/2} ,
\end{eqnarray}
where $\bm{D}$ is the diagonal of $\bm{K}$. The importance of
such scaling and why it helps has been explained for an
idealized one-dimensional problem in \cite{TopKrylov2006}. We
further reduce the condition number of the system matrix and
the number of iterations for convergence by applying an
incomplete Cholesky preconditioner with zero fill-in to the
explicitly rescaled system,
\begin{eqnarray}\nonumber
  \tilde{\bm{K}} & \approx & \tilde{\bm{L}}\tilde{\bm{L}}^T .
\end{eqnarray}

The finite element analysis in topology optimization requires
the solution of a sequence of (usually) symmetric linear
systems. In each optimization step, the algorithm updates the
element densities, and after the first few optimization steps
the changes in the design variables tend to be small from one
optimization step to the next. Hence, the optimization leads to
small changes from one linear system to the next, and the
search space generated for one linear system provides important
information for subsequent linear systems. First, the solution
of one system can be used as an initial guess for the next
system, reducing the initial residual. Second, an approximate
invariant subspace derived from the Krylov space generated for
one linear system can be used for subsequent linear systems,
improving the convergence rate of the iterative solver. This is
the basic idea of Krylov subspace recycling; however, other
subspaces may also be used for 'recycling' \cite{Parks2006}.
Since the linear systems discussed in this paper are symmetric,
we use the Recycling MINRES algorithm (RMINRES) for Krylov
subspace recycling \cite{TopKrylov2006}. Unfortunately, solving
a sequence of problems on meshes that change periodically makes
recycling more difficult. Although it is not hard to map
relatively smooth eigenvectors from a mesh to an updated mesh,
the combination of mesh adaptation and preconditioning seems to
give accuracy problems. Recycling is still effective for AMR;
however, it is not nearly as beneficial as on a static mesh,
and its improvement is work in progress.

\section{Implementation Issues} \label{sec:impls}

For the implementation of adaptive mesh refinement, we use the
\libmesh\  library \cite{libmesh} developed at the University
of Texas at Austin and the Technische Universit\"{a}t
Hamburg-Harburg. The \libmesh\ library provides a C++ framework
for numerical simulations of partial differential equations on
serial and parallel platforms. It supports one-dimensional,
two-dimensional, and three-dimensional finite element and
finite volume simulations on adaptive meshes. The \libmesh\
software uses \petsc\ \cite{petsc-web, petsc-manual} for the
solution of linear systems on both serial and parallel
platforms. However, we use our own custom linear solvers with
Krylov subspace recycling and preconditioners as detailed in
Section~\ref{sec:iter_sol} and references
\cite{TopKrylov2006,Parks2006}. For compatibility with the
\libmesh\ package we have implemented the RMINRES method in the
\petsc\ framework. For the incomplete Cholesky preconditioner
we used routines provided by \petsc.

We have developed two-dimensional and three-dimensional
topology optimization algorithms on top of \libmesh. Currently,
we use element-based design variables, the SIMP method for
material interpolation \cite{Bendsoe1999SIMP}, the OC method
for optimization \cite{BendsoeBk2003, Bendsoe1988}, and
Sigmund's filter technique \cite{SigmundThesis1994,
Sigmund1997, Sigmund1998} with some modifications for dealing
with adaptive refinement. Following Stainko \cite{stainko2006},
we make a small modification in Sigmund's filter for a
nonuniform mesh. Sigmund's filter takes a distance and density
weighted average of the sensitivities of all elements in a
certain radius as
\begin{equation}
\widehat{\frac{\partial c}{\partial \rho_e}} = \frac{1}{\rho_e \sum_d H_{de}} \sum_d \rho_d H_{de} \frac{\partial c}{\partial
  \rho_d},
\label{eq:sig_filter}
\end{equation}
where $\partial c/\partial \rho_e$ is the sensitivity of the
compliance with respect to the density of element $e$, and
$H_{de}$ is a distance weight defined as
\begin{equation}
H_{de} = \max \{ \rmin - \mbox{dist}(d,e), 0 \}.
\end{equation}
The parameter $\rmin$ is a given radius for the filter (for the
work reported here we use $\rmin = \radiusamr$), and
$\mbox{dist}(d,e)$ is the distance between the centers of
elements $d$ and $e$. For a nonuniform mesh, we take the
variation of element size into account by using the element
volume to redefine the weight in the filter  \cite{stainko2006}
as
\begin{equation}
\widehat{\frac{\partial c}{\partial \rho_e}} = \frac{1}{\rho_e \sum_d H_{de} V_d } \sum_d \rho_d H_{de} V_d \frac{\partial
  c}{\partial \rho_d}.      \label{eq:mod_filter}
\end{equation}
The filter radius $\rmin$ is often a length scale independent
of the mesh representation. Notice that the filter will be
effectively deactivated if its size is smaller than that of the
smallest element, i.e., no element has any neighbors within
distance $\rmin$.

\begin{figure}
\begin{center}
\includegraphics[scale=1]{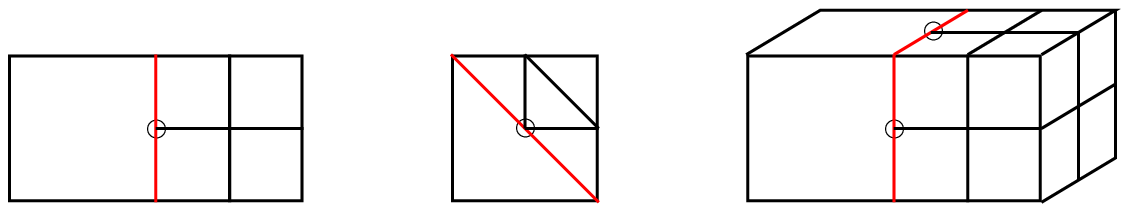} \\
(a) \\[0.2in]
\includegraphics[scale=1]{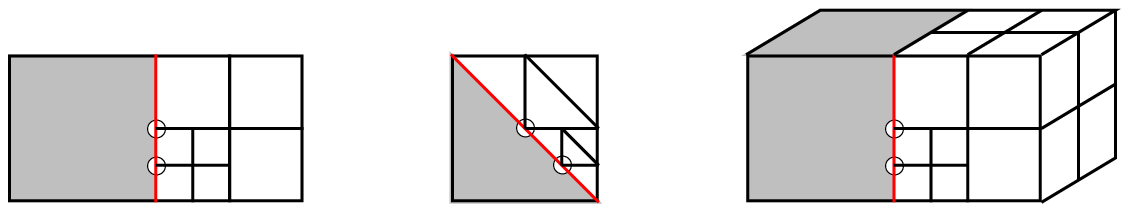} \\
(b)
\end{center}
\caption{Mesh incompatibility with examples of quad, triangle and hex elements: (a) level-one mesh incompatibility marked by
circled nodes; (b) level-two mesh incompatibility marked by circled nodes. We allow level-one
mesh incompatibility (see Section~\ref{sec:impls}), but we avoid level-two and higher incompatibility by refining the gray
coarse elements and by not derefining their children elements if these gray elements result from a potential derefinement.}
\label{fig:mesh_incomp}
\end{figure}

Because of the hierarchical data structure of \libmesh, we must
allow level-one mesh incompatibility. However, we avoid
level-two and higher mesh incompatibility. For example, if the
configuration in Figure~\ref{fig:mesh_incomp}(b) would result
from mesh refinement, we refine the gray elements as well. If
the configuration, in particular the gray elements, in
Figure~\ref{fig:mesh_incomp}(b) would result from mesh
derefinement, we avoid the derefinement. In this way, we limit
mesh incompatibility to level-one mesh incompatibility. As
indicated by the circled nodes in
Figure~\ref{fig:mesh_incomp}(a), level-one mesh incompatibility
results in hanging nodes. The \libmesh\ package handles those
hanging nodes by using the projection method to enforce
constraints in the stiffness matrix. We divide the degrees of
freedom (DOFs) into two groups. Group one consists of all the
unconstrained DOFs, and group two consists of the constrained
DOFs on the hanging nodes. The constrained DOFs can be computed
by linear interpolation from unconstrained DOFs. If we define
vector $\tilde{\bm{u}}$ on the unconstrained DOFs, then
\begin{equation}
\bm{u} = \colvec{\tilde{\bm{u}}}{\bm{P}\tilde{\bm{u}}} =
\matt{\bm{I}}{\bm{0}}{\bm{P}}{\bm{0}} \colvec{\tilde{\bm{u}}}{\bm{0}}
\end{equation}
is the mapping of $\tilde{\bm{u}}$ to all the DOFs, where $\bm{P}$ is the interpolation matrix.
We compute $\hat{\bm u}$ by solving the projected system
\begin{equation}
\matt{\bm{I}}{\bm{P}^T}{\bm{0}}{\bm{0}} \bm{K}
\matt{\bm{I}}{\bm{0}}{\bm{P}}{\bm{0}} \hat{\bm{u}}
= \matt{\bm{I}}{\bm{P}^T}{\bm{0}}{\bm{0}} \bm{f}.
\label{eq:constrained_system}
\end{equation}
Since \libmesh\ does not drop the constrained DOFs in the linear system,
the projected system in (\ref{eq:constrained_system}) is singular
when there is any hanging node. Krylov subspace
methods can handle such singularities as long as the right hand side is consistent, but these singularities may cause problems
for preconditioners. To avoid the singularities, we set the
diagonal entries in the matrix corresponding
to the constrained DOFs to $1$ and solve
\begin{equation}
\left( \matt{\bm{I}}{\bm{P}^T}{\bm{0}}{\bm{0}} \bm{K} \matt{\bm{I}}{\bm{0}}{\bm{P}}{\bm{0}}
+ \matt{\bm{0}}{\bm{0}}{\bm{0}}{\bm{I}}   \right) \hat{\bm{u}}
= \matt{\bm{I}}{\bm{P}^T}{\bm{0}}{\bm{0}} \bm{f}.
\end{equation}
In the end, we recover the constrained DOFs using the interpolation
matrix:
\begin{equation}
\bm{u} = \matt{\bm{I}}{\bm{0}}{\bm{P}}{\bm{0}} \hat{\bm{u}}.
\end{equation}

\section{Results and Discussion} \label{sec:results}

We solve three problems to demonstrate the improvements of our
new AMR scheme and verify that the computed designs using AMR
meshes are equivalent to designs on uniform fine meshes.

For the first (2D) test problem, we compute the optimal design
on a uniform fine mesh and on an adaptively refined mesh with
both our AMR scheme and an approach following references
\cite{Costa2003,stainko2006}. The highest level of refinement
in the AMR meshes has the same element size as that in the
uniform mesh. The results show that our scheme computes a
solution equivalent to the optimal design on the uniform fine
mesh (within a small tolerance), while the alternative AMR
approach from \cite{Costa2003,stainko2006} does not. Moreover,
the experiments elucidate how this suboptimal design arises
from the strategy to only refine the results from a fixed
coarse mesh.

For the second (3D) test problem, we compare the optimal design
using an adaptive mesh and our AMR strategy with the optimal
design on a uniform fine mesh for a relatively simple
cantilever beam problem. Again the maximum refinement level for
the AMR mesh has the same element size as the uniform, fine
mesh. We also use this test to show that our AMR strategy leads
to faster convergence in both the linear iterations (finite
element solution) and the nonlinear iterations (design
optimization) as well as to a significant reduction in runtime
(about a factor of three).

In the third test problem (also 3D), we compare the optimal
design using an adaptive mesh and our AMR strategy with the
optimal design on a uniform, fine mesh for a more complicated
design problem. For all three test problems our AMR strategy
leads to essentially the same design as is obtained on a
uniform mesh, but at significantly reduced cost.

To evaluate the relative difference between two designs, we
need a quantitative measure. We define the the relative
difference between two designs as
\begin{equation}
  D(\rho^{(1)},\rho^{(2)}) = \frac{ \int_{\Omega} | \rho^{(1)}-\rho^{(2)} | d\Omega }{ \int_{\Omega} \rho^{(1)} d\Omega }.
\label{eq:diff_measure}
\end{equation}
We take $\rho{(1)}$ to indicate the design on the uniform fine
mesh, and $\rho{(2)}$ to indicate the design on the AMR mesh.
This difference can be computed in a number of ways. To
simplify comparison, we take the uniform mesh to be the AMR
mesh with maximum refinement at every point. So, we refine the
AMR mesh to the same fine mesh level at every point (without
changing the design), and then evaluate the 1-norm of the
difference between the designs.

\subsection*{Test 1: 2D cantilever beam}
\begin{figure}
\begin{center}
\includegraphics{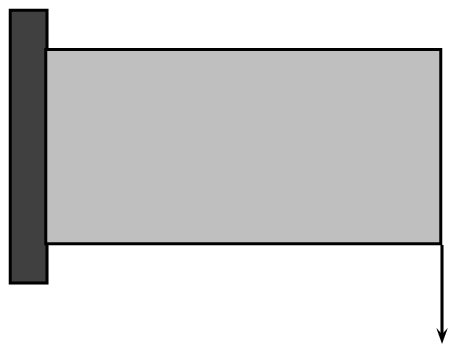}\\
(a)\\[0.2in]
\includegraphics[scale=0.3]{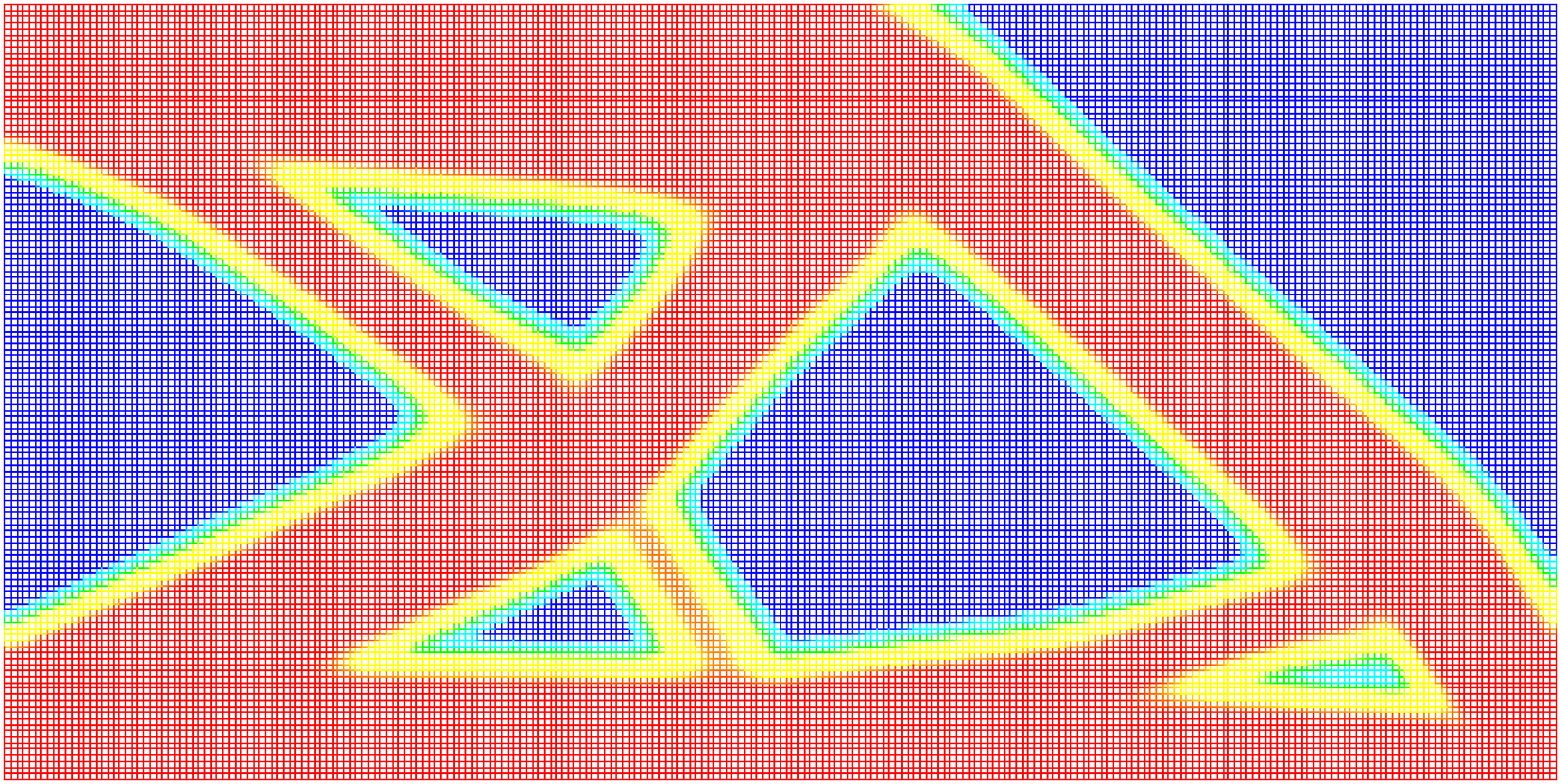} \includegraphics[scale=0.5]{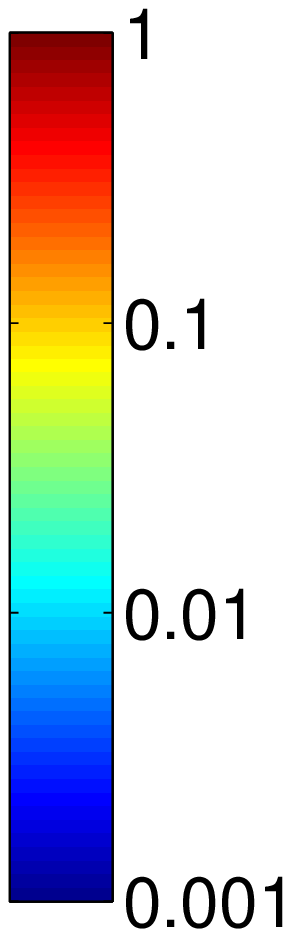} \\
(b)\\[0.2in]
\includegraphics[scale=0.3]{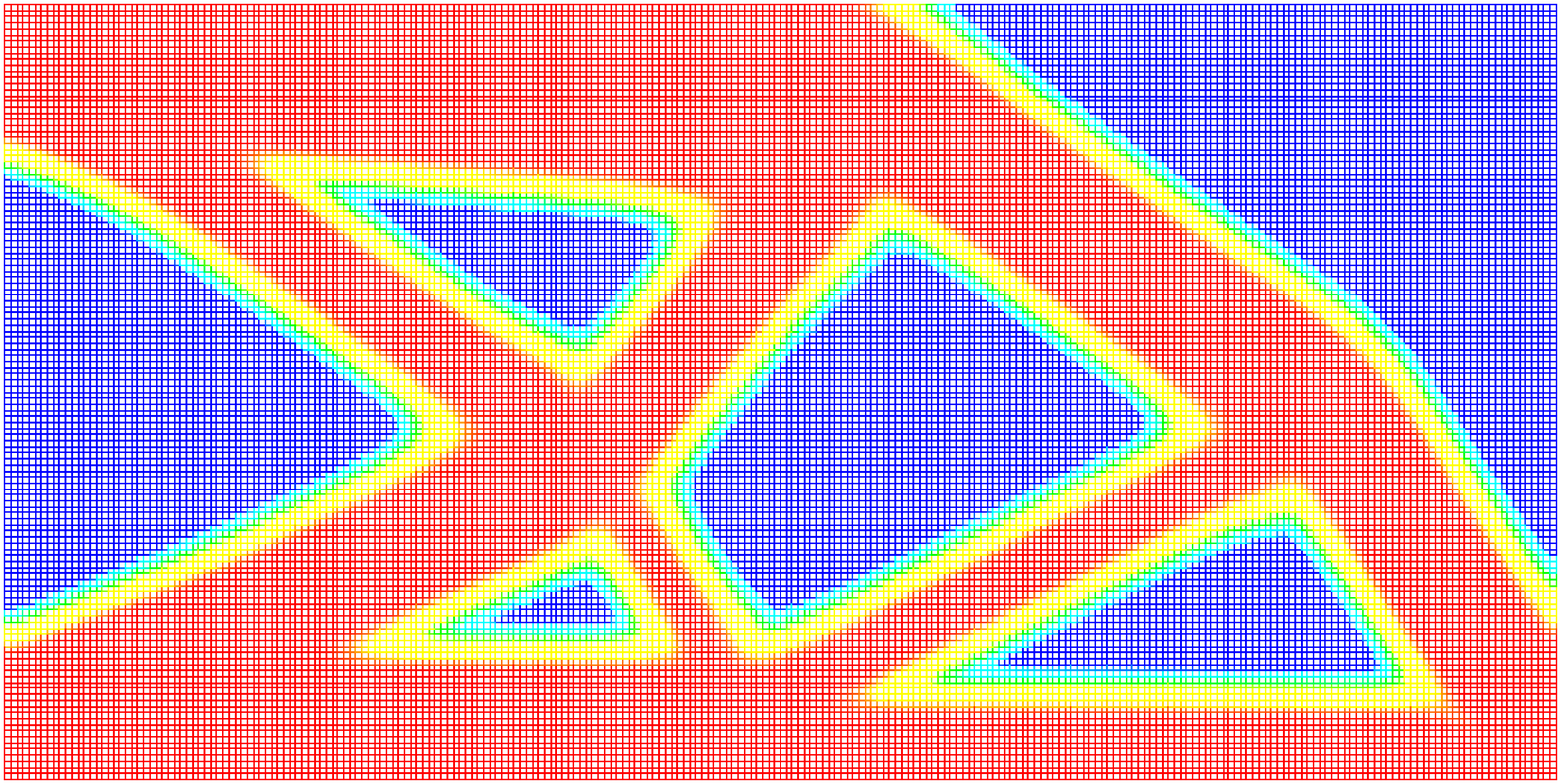} \includegraphics[scale=0.5]{figures/colormap/colormap_log} \\
(c)
\end{center}
\caption{Topology optimization on a $256\times128$ \emph{uniform mesh}: (a) problem configuration (the volume constraint $V_0$ is $50\%$ of the domain volume);
  (b) an intermediate design; (c) final converged design.}
\label{fig:noamr}
\end{figure}

We compute the optimal design for the 2D beam problem shown in
Figure \ref{fig:noamr}(a). We first compute the design on a
uniform, fine mesh. Figure \ref{fig:noamr}(b) shows an
intermediate result, and Figure \ref{fig:noamr}(c) the
converged design. Note how the truss member at the lower-right
corner has risen up noticeably from the intermediate result to
the final design. An effective AMR procedure must be able to
capture such an evolution in the design.

\begin{figure}
\begin{center}
\includegraphics[scale=0.25]{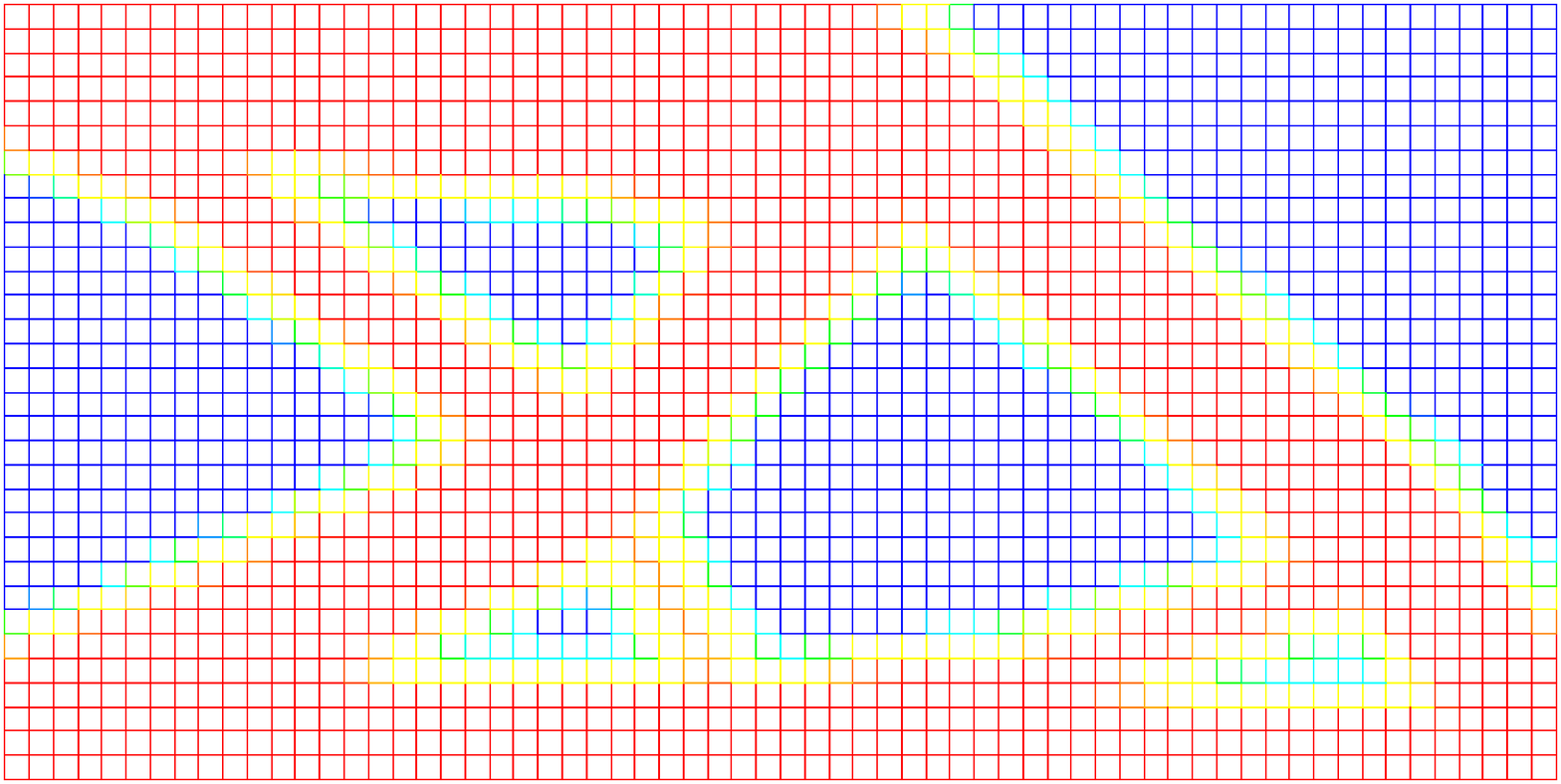} \includegraphics[scale=0.4]{figures/colormap/colormap_log} \\
(a)\\[0.2in]
\includegraphics[scale=0.25]{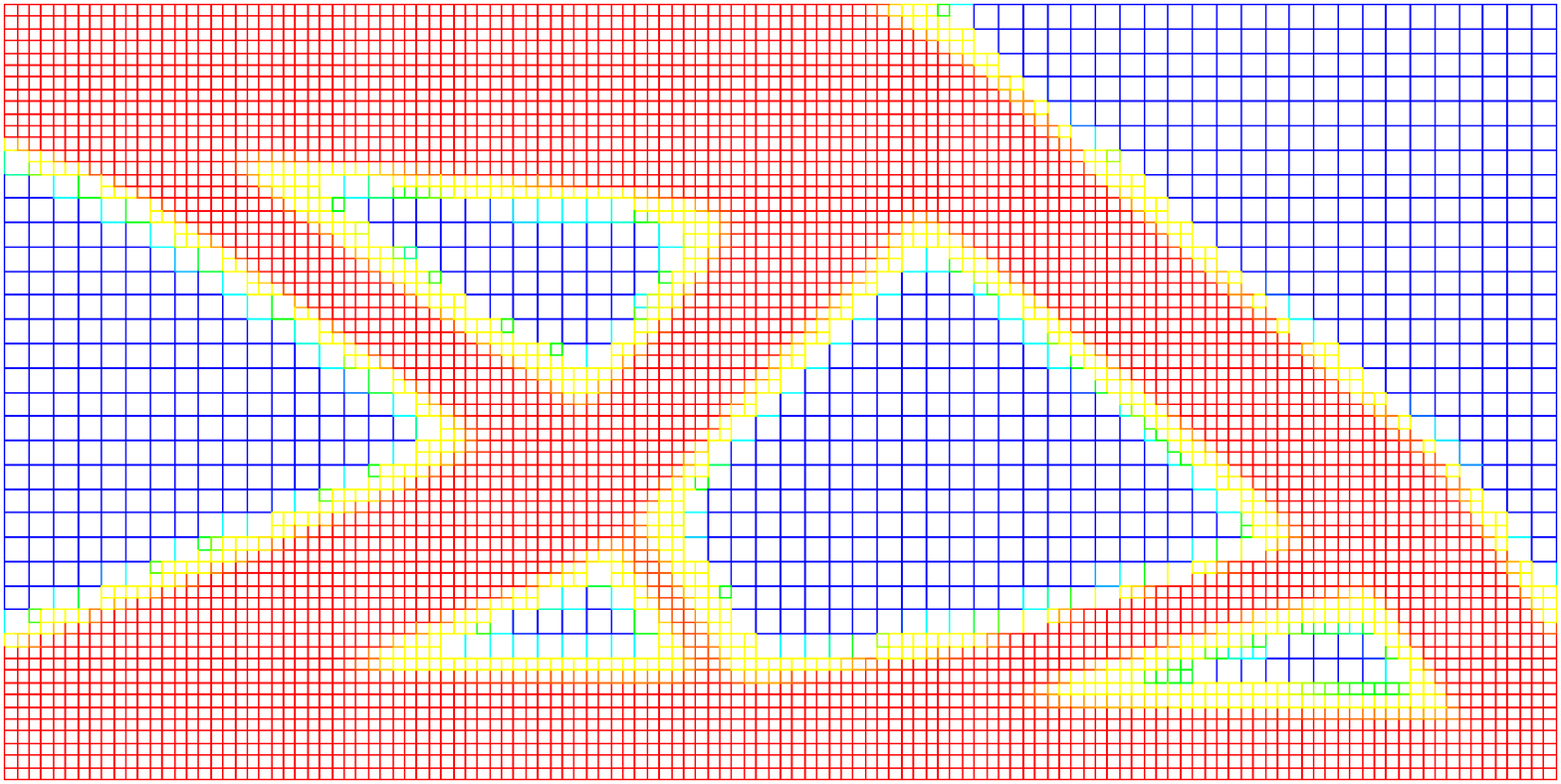} \includegraphics[scale=0.4]{figures/colormap/colormap_log} \\
(b)\\[0.2in]
\includegraphics[scale=0.25]{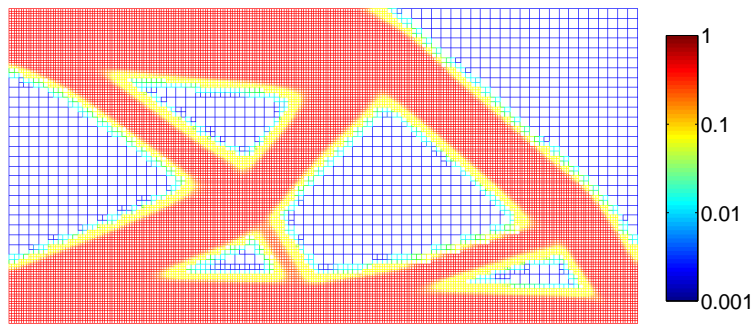} \includegraphics[scale=0.4]{figures/colormap/colormap_log} \\
(c)
\end{center}
\caption{Topology optimization on \emph{an adaptive mesh with only refinement on each level}: (a) converged result on the coarsest mesh with 2048 elements;
  (b) converged result on the intermediate mesh with 5675 elements;
  (c) converged result on the final mesh with 20216 elements. Note the
  undesirable position of the truss member near the lower right corner, which remains nearly invariant during the evolution of the design.}
\label{fig:single_refine}
\end{figure}

Next, we solve the same problem following the strategy
mentioned in references \cite{Costa2003,stainko2006}. We start
with a relatively coarse mesh ($64 \times 32$), and obtain the
converged solution to the topology optimization problem shown
in Figure \ref{fig:single_refine}(a). Then, we refine the mesh
according to this coarse level result and we solve the
optimization problem on this locally refined mesh until
convergence, obtaining the solution shown in Figure
\ref{fig:single_refine}(b). Next, we refine the mesh and solve
again. Finally, we obtain the result on the finest mesh shown
in Figure \ref{fig:single_refine}(c). The truss member at the
lower-right corner has remained roughly at its original
position on the coarsest mesh in spite of the high resolution
of the design, causing the resulting design to differ
significantly from the optimal design obtained on the uniform,
fine mesh, even though the smallest element size in the meshes
in Figure \ref{fig:single_refine}(c) is the same as the element
size for the uniform mesh shown in Figure \ref{fig:noamr}. This
difference in material distribution is caused by the fine mesh
discretization being confined to the region identified by the
coarse mesh design. The mesh adaptation strategy only allows
the fine mesh computation to refine the coarse mesh design. It
does not allow the fine mesh computation to alter the design
substantially, even if more accurate fine mesh computations
indicate a better design. An additional problem is that the
initial mesh needs to be relatively fine, such as the one in
Figure \ref{fig:single_refine}(a), because a coarser initial
mesh would lead to very poor solutions as the filter would be
inactive at that mesh resolution. In this case, mesh adaptation
with refinement only leaves fine elements that could be
derefined for efficiency in the void regions.

\begin{figure}
\begin{center}
\includegraphics[scale=0.25]{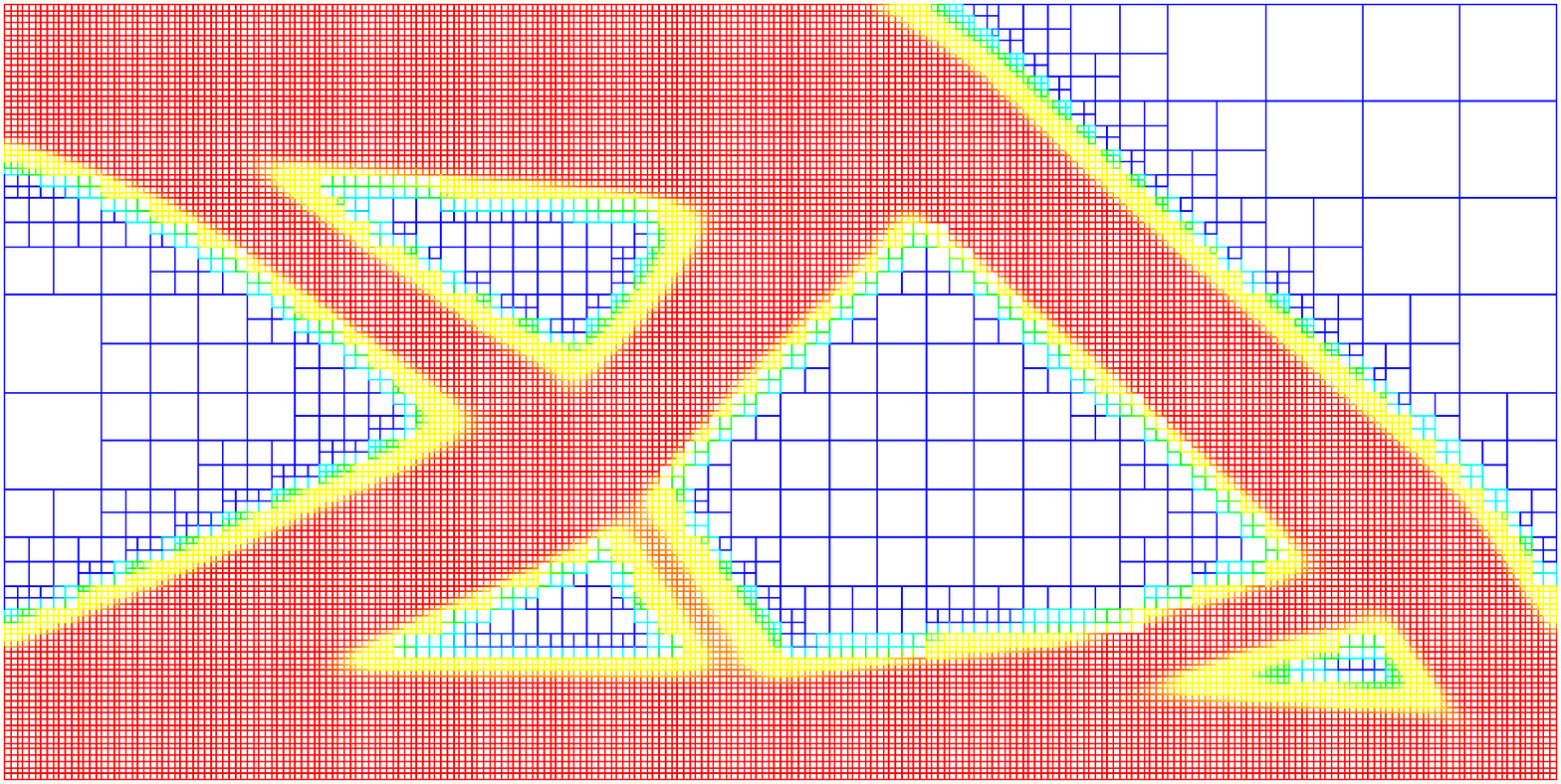} \includegraphics[scale=0.4]{figures/colormap/colormap_log} \\
(a)\\[0.2in]
\includegraphics[scale=0.25]{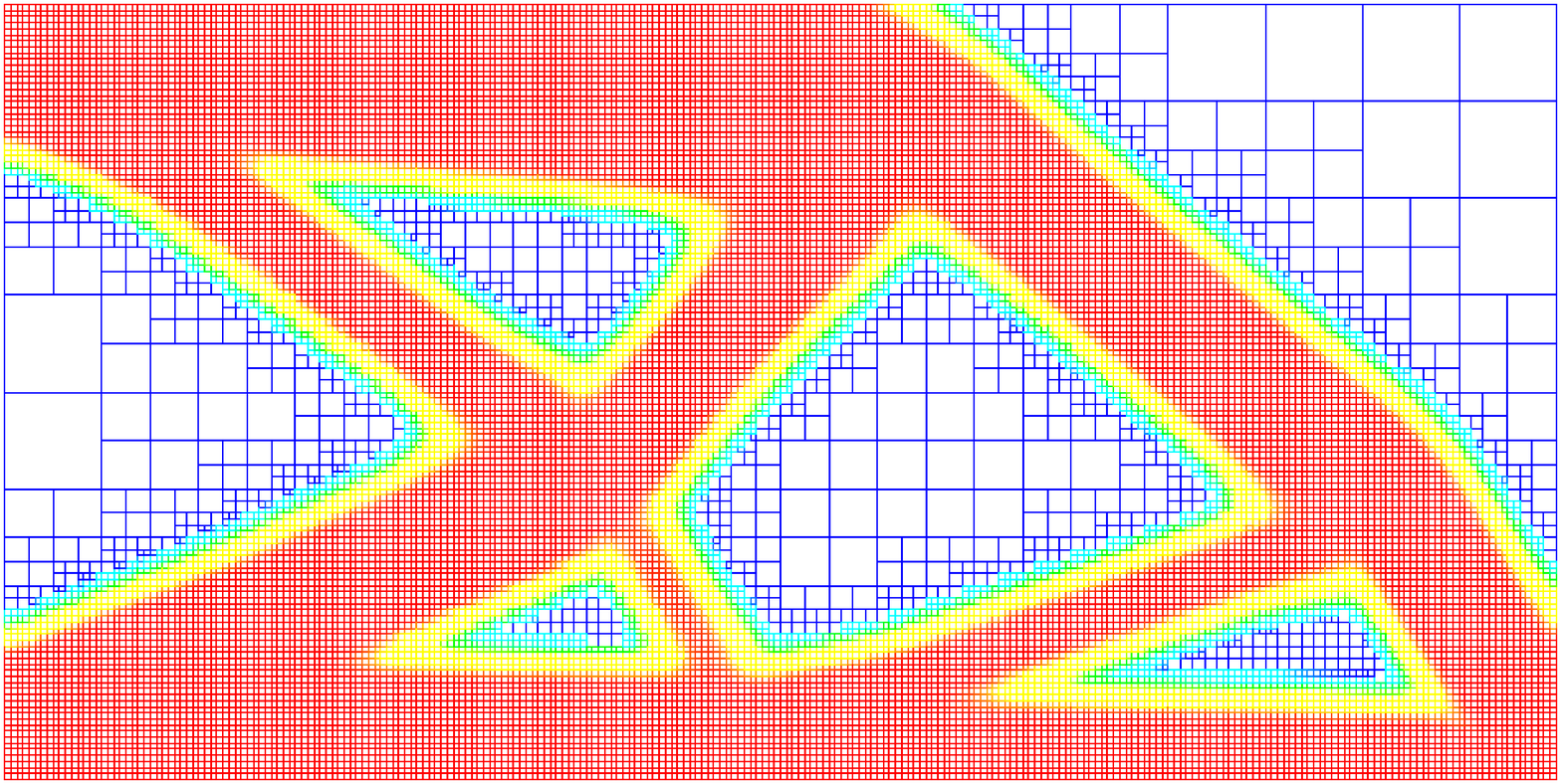} \includegraphics[scale=0.4]{figures/colormap/colormap_log} \\
(b)\\[0.2in]
\includegraphics[scale=0.25]{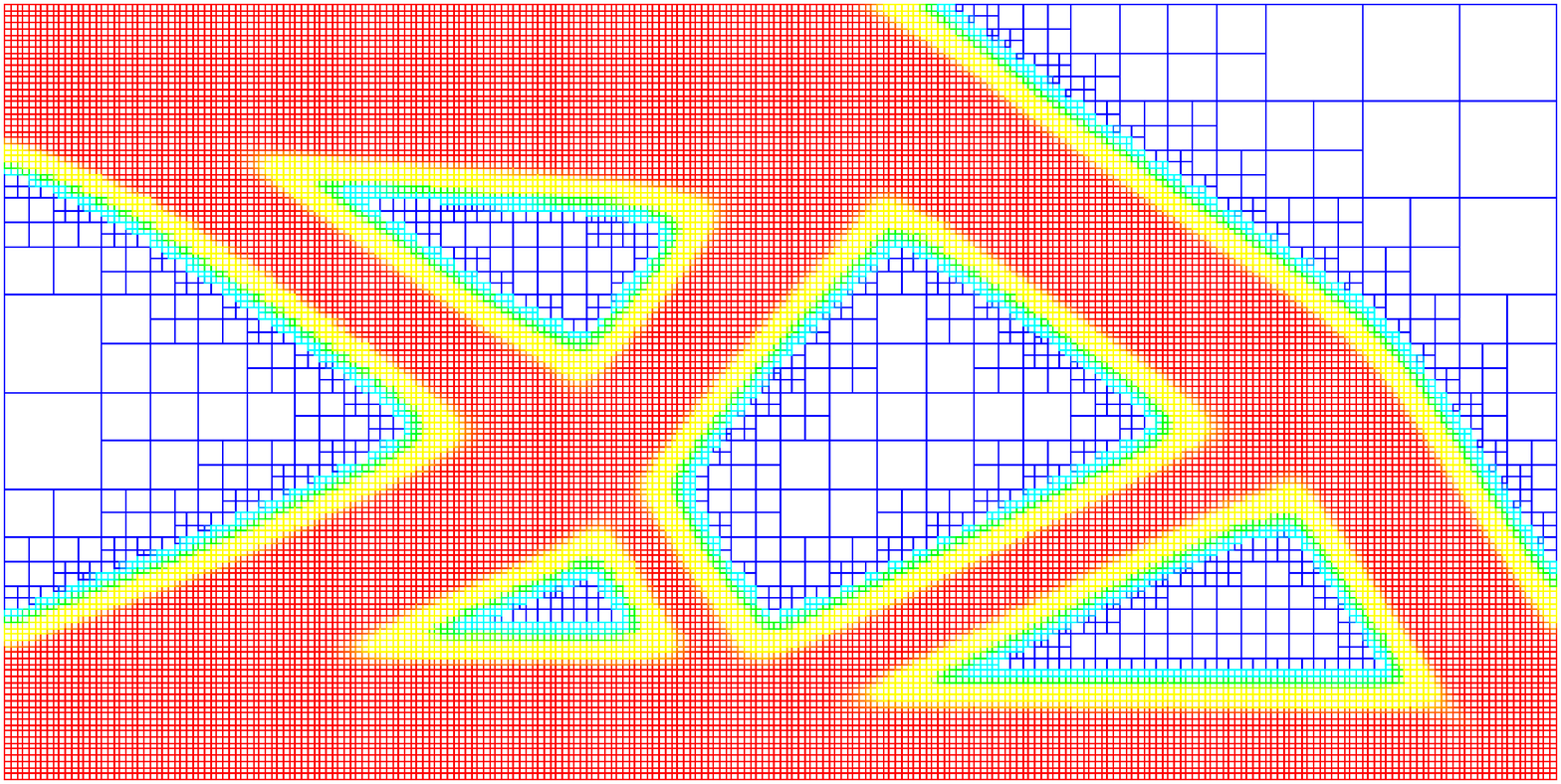} \includegraphics[scale=0.4]{figures/colormap/colormap_log} \\
(c)
\end{center}
\caption{Topology optimization on \emph{an adaptive mesh with continual
dynamic refinement and derefinement on each level}: (a)--(b) intermediate designs; (c) final converged design on a nonuniform mesh with
   $25229$ elements, whose finest resolution is the same as the resolution of the uniform mesh in Figure \ref{fig:noamr}. Notice that the truss member of the lower-right corner moves up as the AMR procedure progresses.}
\label{fig:multi_refine}
\end{figure}

\begin{table}
\caption{\label{tab:AMRCosta} Mesh adaptation scheme following Costa and Alves \cite{Costa2003}.}
\begin{center}
\begin{tabulary}{\textwidth}{CCCCC}
\toprule
opt. step & $\ell_{\max}$ & \#elem & \#unknowns & relative $L_1$ difference \\
\midrule
80 & 1 & 2048 & 4290 & 21.16\% \\
125 & 2 & 5675 & 11948 & 19.53\% \\
200 & 3 & 20216 & 41824 & 19.42\% \\
300 & 3 & 20216 & 41824 & 19.60\% \\
400 & 3 & 20216 & 41824 & 19.66\% \\
500 & 3 & 20216 & 41824 & 19.66\% \\
600 & 3 & 20216 & 41824 & 19.65\% \\
700 & 3 & 20216 & 41824 & 19.67\% \\
800 & 3 & 20216 & 41824 & 19.66\% \\
900 & 3 & 20216 & 41824 & 19.65\% \\
1000 & 3 & 20216 & 41824 & 19.67\% \\
\bottomrule
\end{tabulary}
\end{center}
\end{table}

\begin{table}
\caption{\label{tab:AMRSEG} Dynamic AMR scheme}
\begin{center}
\begin{tabulary}{\textwidth}{CCCCC}
\toprule
opt. step & $\ell_{\max}$ & \#elem & \#unknowns & relative $L_1$ difference \\
\midrule
27 & 1 & 2048 & 4290 & 21.13\% \\
37 & 2 & 6113 & 12786 & 19.80\% \\
100 & 3 & 24707 & 50560 & 17.77\% \\
200 & 3 & 25040 & 51242 & 13.00\% \\
300 & 3 & 25184 & 51550 & 5.00\% \\
367 & 3 & 25229 & 51654 & 0.17\% \\
\bottomrule
\end{tabulary}
\end{center}
\end{table}

Now, we solve the same problem, starting with the coarse mesh
of Figure~\ref{fig:single_refine}(a), but following our AMR
strategy (Section \ref{sec:dynamic_amr}). We allow multiple
mesh adaptations on any level and we maintain a layer of fine
elements on the void side of the solid/void interface. This
leads to the results shown in Figure~\ref{fig:multi_refine},
with two intermediate results shown in
Figures~\ref{fig:multi_refine}(a) and (b), and the final
converged result in Figure~\ref{fig:multi_refine}(c). Note how
the truss member at the lower-right corner moves up as the
optimization progresses, just as for the evolution of
intermediate designs on the uniform mesh (reference). The
figures also demonstrate how the mesh changes smoothly with the
changes in material distribution. The smallest element size in
the AMR meshes in Figure \ref{fig:multi_refine} is the same as
the element size for the uniform mesh shown in
Figure~\ref{fig:noamr} (reference). Compared with the final
solution shown in Figure~\ref{fig:single_refine}, the solution
obtained with our AMR strategy is closer to the solution
obtained on the uniform mesh (reference). Indeed, based on the
metric in (\ref{eq:diff_measure}), the relative difference
between the designs in Figure \ref{fig:single_refine}(c) and
Figure \ref{fig:noamr}(c) is $19.6 \%$, while the relative
difference between the designs in
Figure~\ref{fig:multi_refine}(c) and Figure \ref{fig:noamr}(c)
is only $0.168\%$. Furthermore, derefinement results in coarser
elements in the void regions, cf.
Figures~\ref{fig:single_refine}(c) and
\ref{fig:multi_refine}(c). These results are summarized in
Tables~\ref{tab:AMRCosta} and \ref{tab:AMRSEG}, where
$\ell_{\max}$ refers to the highest refinement level present in
the mesh.

\subsection*{Test 2: 3D cantilever beam}
\begin{figure}
\begin{center}
\includegraphics[scale=0.8]{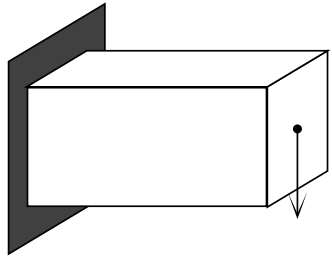}
\end{center}
\caption{3D cantilever beam example with domain scale 2:1:1.}
\label{fig:cantilever3d_config}
\end{figure}

\begin{figure}
\begin{center}
\includegraphics[scale=0.25]{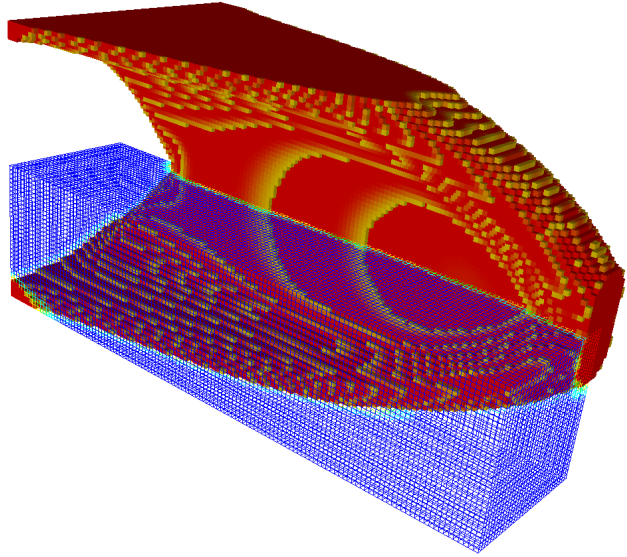} \\
(a) \\
\includegraphics[scale=0.25]{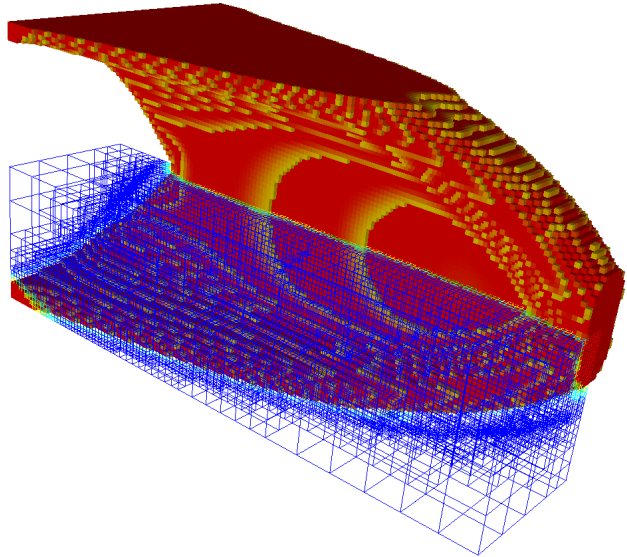} \\
(b)
\end{center}
\caption{Final solutions of the 3D cantilever beam problem
(Figure~\ref{fig:cantilever3d_config}) obtained using symmetry
on a quarter of the domain as indicated by the mesh: (a) final
solution on a fixed uniform mesh with $128\times32\times32$
elements; (b) AMR solution on a mesh with $57173$ elements; the
finest local resolution is the same as that of the uniform,
fine mesh.} \label{fig:cantilever3d_amr}
\end{figure}

\begin{figure}
\begin{center}
\begin{tabular}{cc}
\multicolumn{2}{c}{\includegraphics[scale=0.5]{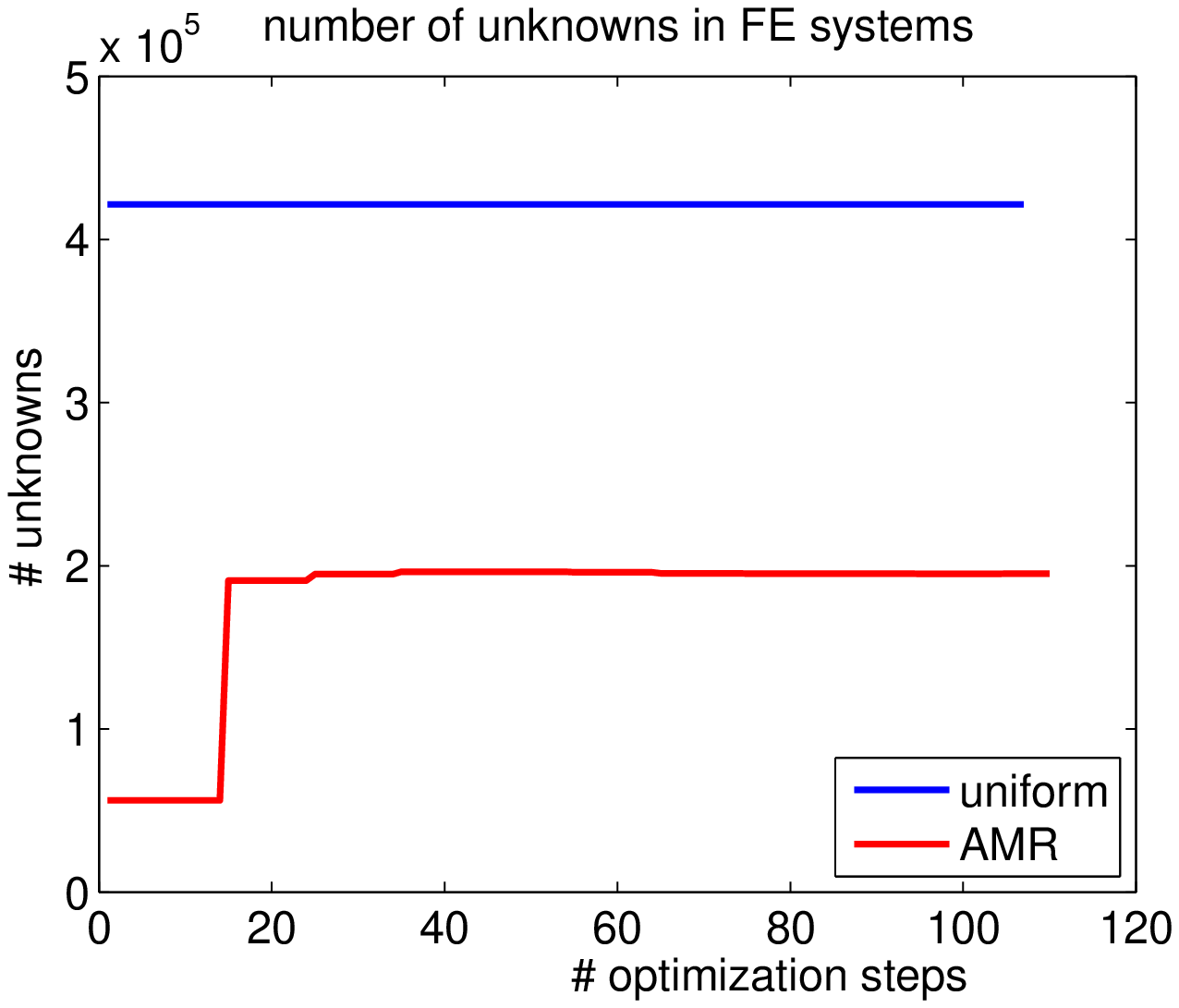}} \\
\multicolumn{2}{c}{(a)} \\
\includegraphics[scale=0.5]{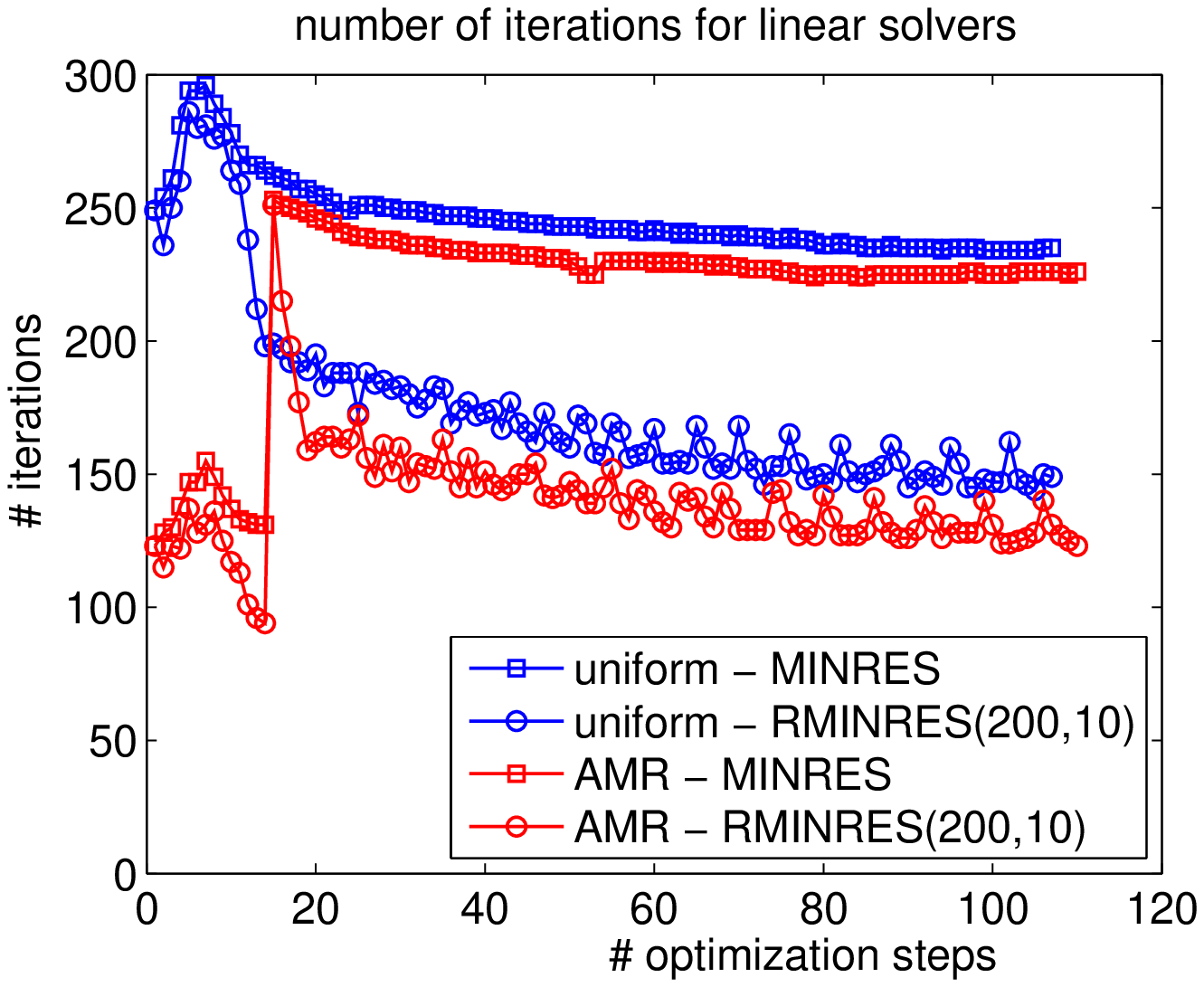} &
\includegraphics[scale=0.5]{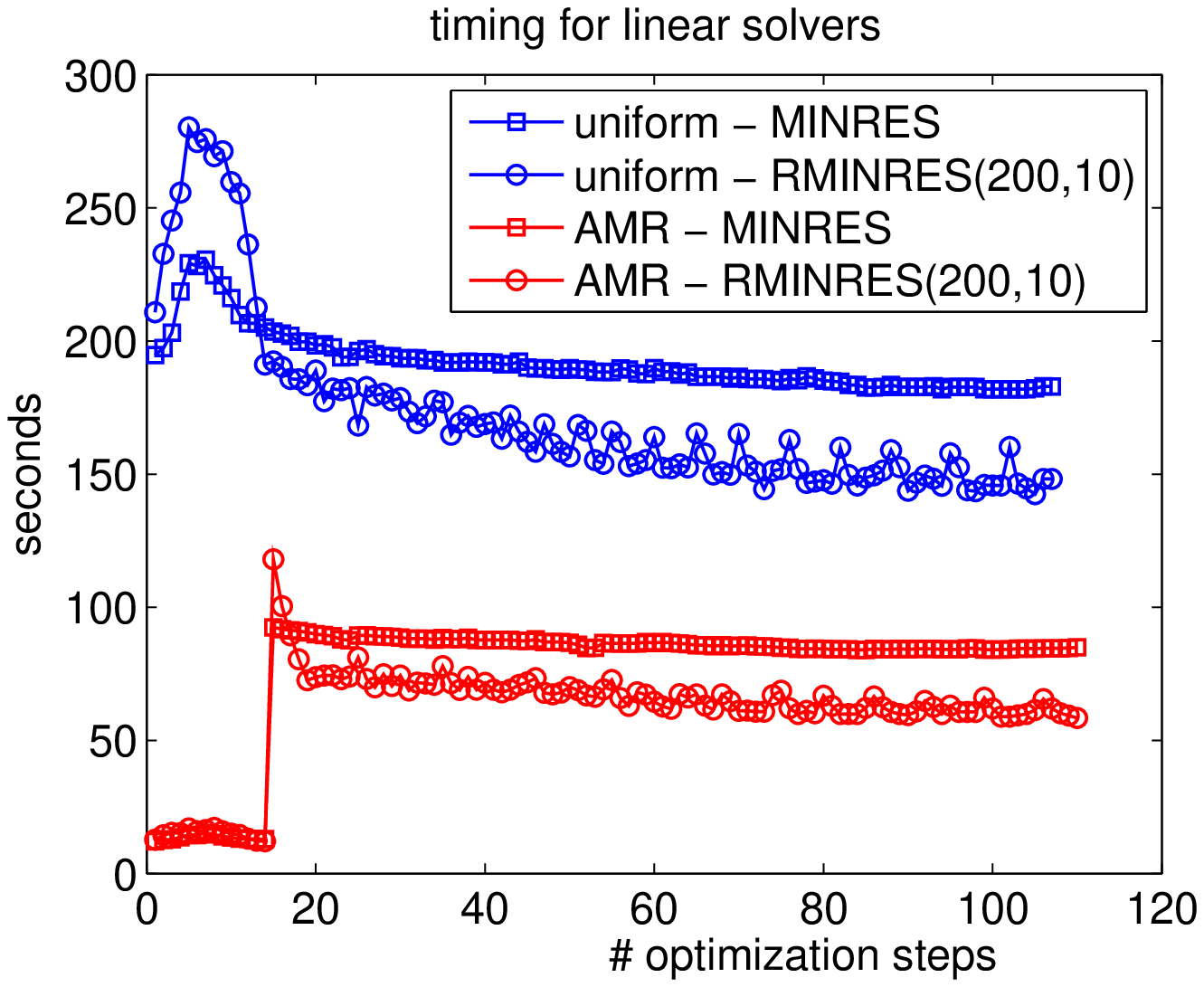} \\
(b) & (c)
\end{tabular}
\end{center}
\caption{Comparison of linear solver statistics for the
cantilever beam design problem on a uniform, fine mesh and on
an adaptive mesh: (a) the number of unknowns in the linear
systems arising from the finite element discretization; (b) the
number of preconditioned MINRES and RMINRES(200,10) iterations
(see below) for each optimization step; (c) solution times with
MINRES and RMINRES(200,10) for the linear systems arising from
finite element discretization at each optimization step. The
parameters $m$ and $k$ in RMINRES(m,k) have the following
meaning. The method recycles an approximate invariant subspace
associated with the smallest $k$ eigenvalues from one linear
system to the next. In the solution of single linear system,
the approximate invariant subspace is updated every $m$
iterations.} \label{fig:amr3d_stats}
\end{figure}

We compute the optimal design for the three-dimensional
cantilever beam, shown in Figure~\ref{fig:cantilever3d_config},
with a volume constraint of $25\%$. Exploiting symmetry, we
discretize only a quarter of the domain. We solve this problem
on a (fixed) uniform mesh with $128\times32\times32$ B8
elements and also following our AMR strategy. The initial mesh
for the AMR-based design has $64\times16\times16$ B8 elements.
The final results are shown in
Figure~\ref{fig:cantilever3d_amr} with
Figure~\ref{fig:cantilever3d_amr}(a) displaying the solution on
a uniform, fine mesh, and Figure~\ref{fig:cantilever3d_amr}(8)
displaying the AMR solution; note the large blocks in parts of
the void region in Figure~\ref{fig:cantilever3d_amr}(b). The
relative difference between these two designs is only
$0.0909\%$ (Eq. (\ref{eq:diff_measure})). We use the
preconditioned, recycling minimum residual solver (RMINRES)
proposed in \cite{TopKrylov2006} to solve the linear systems
arising from the finite element discretization for a given
material distribution. The dimensions of the linear systems of
equations for the adaptive mesh are less than half of those for
the uniform, fine mesh. The difference is even larger early in
the optimization iteration. Moreover, the number of RMINRES
iterations for the linear systems derived from the adaptive
mesh are slightly smaller than those for the uniform, fine mesh
(much smaller early in the optimization), because the adaptive
meshes tend to lead to better conditioned linear systems.
Therefore, using AMR reduces the solution time roughly by a
factor of three; see the statistics in
Figure~\ref{fig:amr3d_stats}.

\begin{figure}
\begin{center}
\includegraphics[scale=0.7]{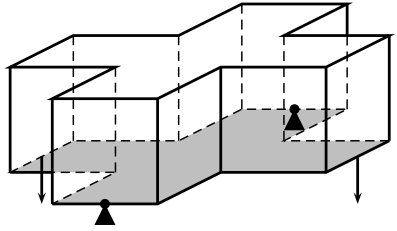}
\end{center}
\caption{A 3D compliance minimization problem in a cross-shaped
domain with the bottom (shaded) front and back ends fixed and
the bottom left and right ends pulled down. The volume
constraint $V_0$ is $20\%$ of the domain volume.}
\label{fig:cross}
\end{figure}

\begin{figure}
\begin{center}
\includegraphics[scale=0.85]{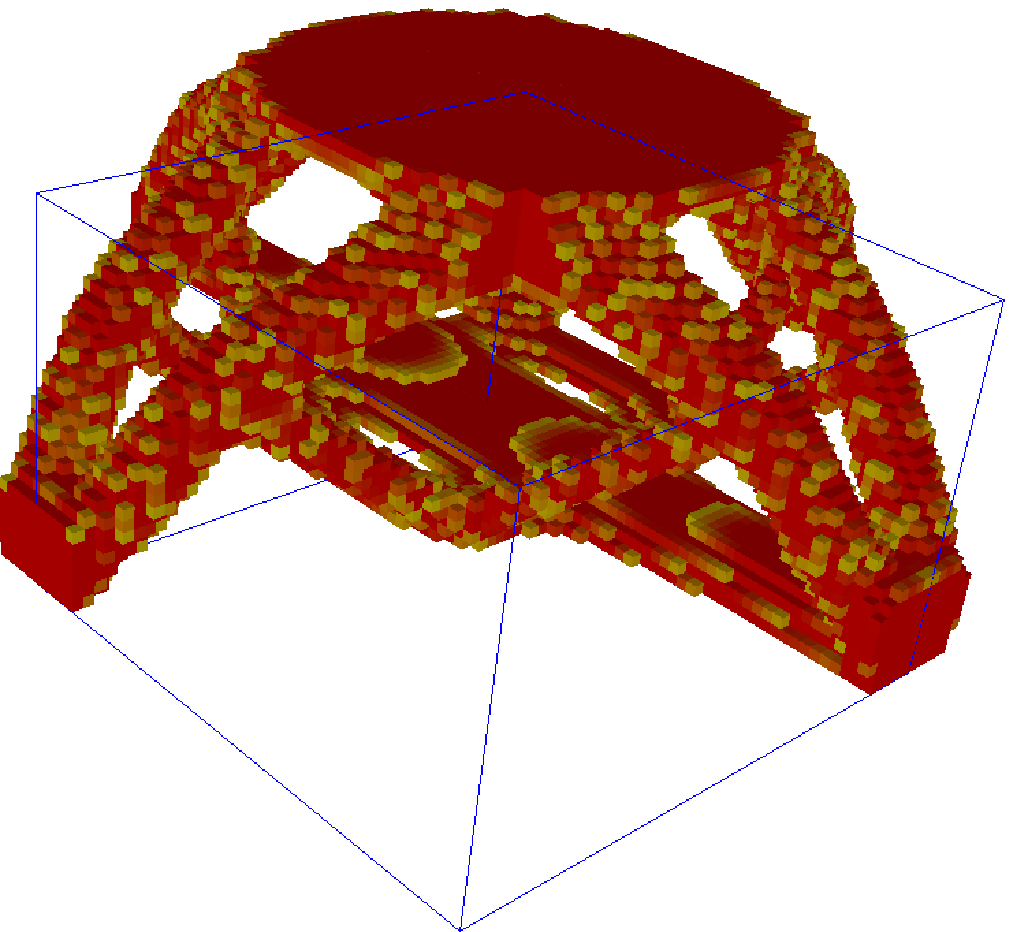} \\[0.3in]
\includegraphics[scale=0.85]{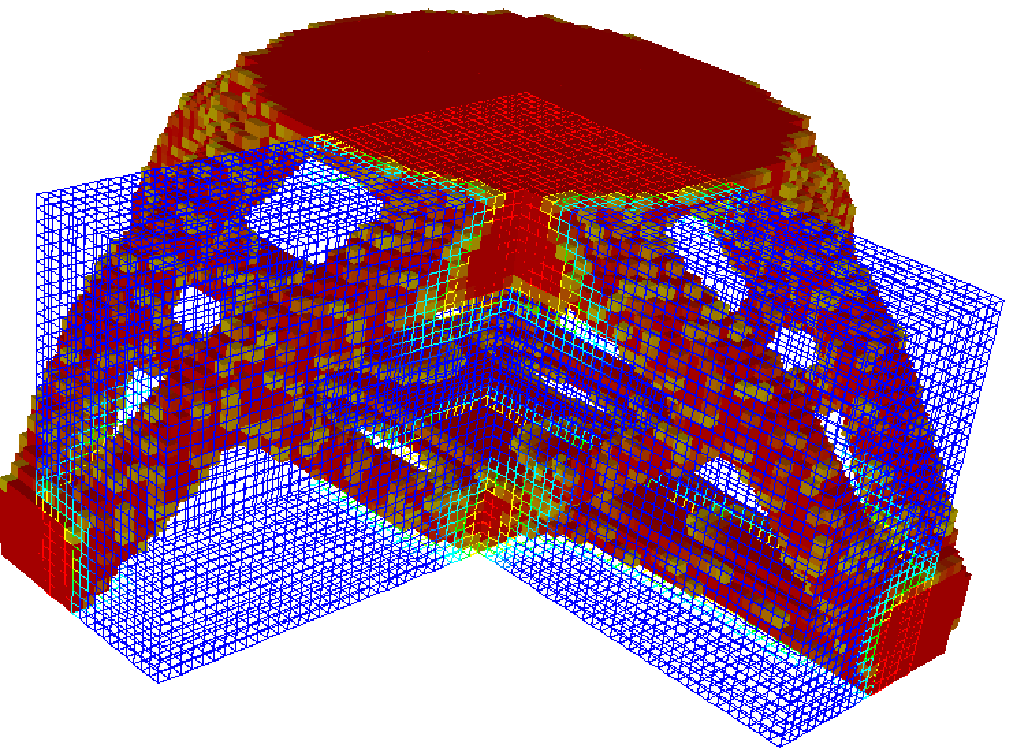}
\end{center}
\caption{The optimal solution to the design problem shown in
Figure \ref{fig:cross} on a uniform finite element mesh with
40960 B8 elements. The quarter-mesh discretization is shown on
the bottom figure.} \label{fig:cross_uniform}
\end{figure}

\begin{figure}
\begin{center}
\includegraphics[scale=0.85]{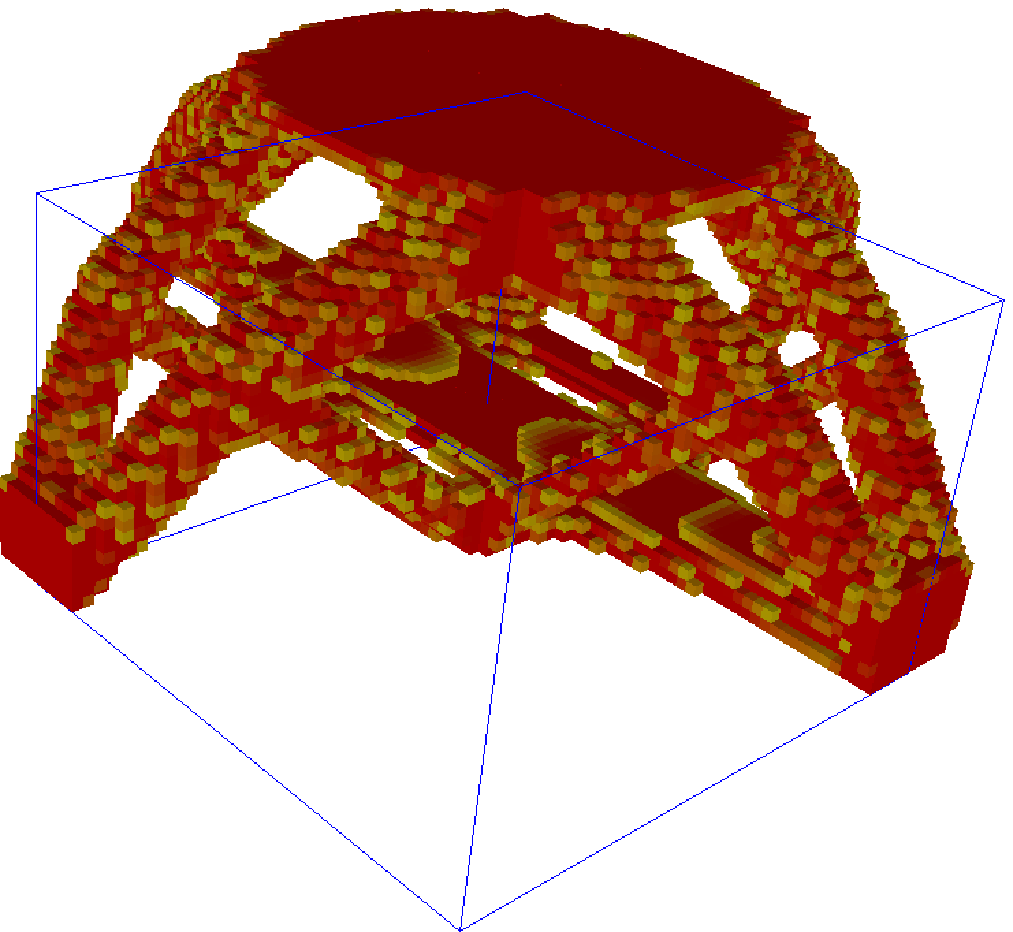} \\[0.2in]
\includegraphics[scale=0.85]{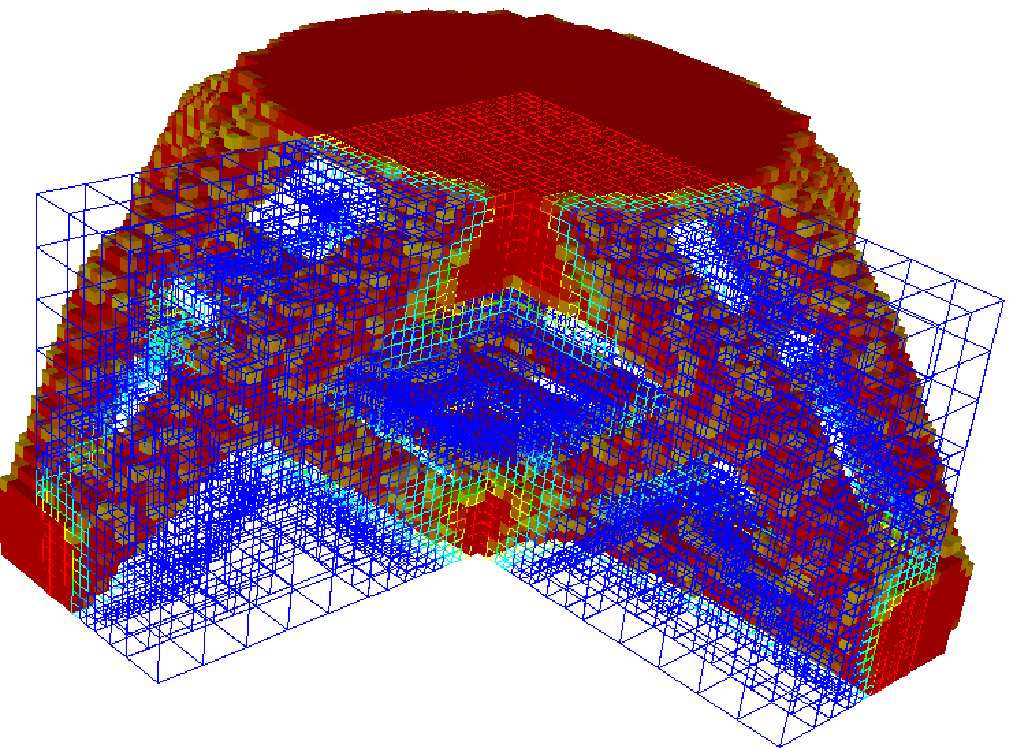}
\end{center}
\caption{The optimal solution to the design problem shown in
Figure \ref{fig:cross} on an adaptively refined mesh. The final
mesh consists of 19736 B8 elements. The quarter-mesh
discretization is shown on the bottom figure.}
\label{fig:cross_amr}
\end{figure}

\subsection*{Test 3: Cross-shaped domain}
We compute the optimal design for the more complex
three-dimensional test problem shown in Figure \ref{fig:cross}.
For the cross-shaped domain, we compute the optimal design
subject to the fixed boundary on the bottom front and back
ends, and two loads on the left and right sides. The maximum
volume allowed is $20\%$ of the domain volume. We solve this
problem both on a uniform mesh and on an adaptive mesh
following our AMR strategy. The results are shown in Figures
\ref{fig:cross_uniform} and \ref{fig:cross_amr}, respectively.
The uniform mesh consists of $40960$ B8 elements, while the
final adaptive mesh consists of only $19736$ B8 elements.
Moreover, the optimization converges in over $200$ steps on the
uniform mesh, but in only $106$ optimization steps on the
adaptive mesh. The adaptive mesh refinement reduces the total
solution time by more than a factor of three to about $30\%$ of
the solution time for the uniform mesh. Nonetheless, the
relative difference between these two designs is only $2.58\%$
(Eq. (\ref{eq:diff_measure})).

\section{Conclusions} \label{sec:conc}

In order to reduce the high computational cost of accurate
three-dimensional designs by topology optimization we use
adaptive mesh refinement. We propose several critical
improvements to the approaches proposed by Costa and Alves
\cite{Costa2003} and Stainko \cite{stainko2006} in order to
attain better designs. In particular, we want to obtain the
same optimal designs that would be obtained on a uniform, fine
mesh with AMR discretization having significantly fewer
elements but the same fine mesh resolution. The purpose of AMR
is to reduce the cost for the (same) optimal design; we do not
want to reduce the quality of designs. For large, complex,
three-dimensional design problems we could not possibly use a
uniform fine mesh at the desired resolution. Our approach
requires a dynamic meshing strategy that involves continual
refinements and derefinements following the strategy laid out
in Section~\ref{sec:dynamic_amr}. Derefinements should also
lead to further efficiency improvement by reducing the number
of elements in void regions, especially for three-dimensional
problems. Using three test problems, we demonstrate that our
AMR algorithm achieves the desired designs that are within a
small tolerance of those obtained on a uniform, fine mesh with
the same finest level. Our AMR strategy significantly reduces
the total runtime, nonlinear, and linear iterations with
respect to using uniform meshes.

Important future work includes error estimation in the finite
element analysis and mesh refinement and derefinement governed
by both considerations of accurate design and error estimation.
In addition, we plan to work on preconditioners that can be
adapted with the mesh (rather than recomputed) and to improve
the convergence rate of Krylov methods with subspace recycling
\cite{TopKrylov2006,Parks2006}. We also intend to extend the
present AMR technique to multiphysics problems
\cite{Carbo2007}.

\section{Acknowledgements}
We are indebted to Hong Zhang and Mat Knepley from Argonne
National Laboratory for their help with the \petsc\ library and
to Roy Stogner, John Peterson, and Benjamin Kirk from the
University of Texas at Austin for their help with the
\libmesh\/ library. We also thank Cameron Talischi for
insigthful discussions which contributed to improve the present
manuscript.

\bibliographystyle{mystyle}
\bibliography{top_amr}

\end{document}